\documentclass{amsart}%
\usepackage{amssymb}
\usepackage{amsmath}
\usepackage{amsfonts}%
\usepackage{graphicx}

	\def\r{\rho}

	\def\<{\langle}
	\def\>{\rangle}
	\def\a{\alpha}
	
	\def\eps{\varepsilon}
	
	\newcommand{\SE}{\setcounter{equation}{0}}
\newtheorem{theorem}{Theorem}
\theoremstyle{plain}

\newtheorem{definition}{Definition}

\newtheorem{lemma}{Lemma}

\newtheorem{proposition}{Proposition}

\numberwithin{equation}{section}
	
\begin{document}

\title[Porous elastic system with second spectrum]{Energy decay analysis for Porous elastic system with microtemperature : A second spectrum approach}
\author{Hamza Zougheib}
\author{ Toufic El Arwadi  }
\author{ Mohammad El-Hindi  }

\address{Hamza Zougheib\newline  Department of Mathematics and Computer
		Science, Faculty of Science, Beirut Arab University, Beirut, Lebanon  } \email{hmz232@student.bau.edu.lb}
		\address{Toufic El Arwadi \newline  Department of Mathematics and Computer
		Science, Faculty of Science, Beirut Arab University, Beirut, Lebanon  } \email{t.elarwadi@bau.edu.lb}
			\address{Mohammad El-Hindi \newline  Department of Mathematics and Computer
		Science, Faculty of Science, Beirut Arab University, Beirut, Lebanon  } \email{mohammadyhindi98@gmail.com}
    
\maketitle

\begin{abstract}
  In this work, we analyze porous elastic system with microtemperature from second spectrum viewpoint. Indeed,  by using the classical Faedo-Galerkin method combined with the a priori estimates, we prove the existence and uniqueness of
a global solution of this problem. Then we prove that this solution is exponentially stable without assuming the condition
of equal wave speeds. Then, we introduce a finite element approximation and we prove that the associated discrete energy decays. Finally, we obtain some a priori error estimates assuming additional regularity on the solution and we present some numerical results
which demonstrate the accuracy of the approximation and the behaviour of the solution
\end{abstract} 

\section{Introduction}
In later a long time, there are so numerous mathematical researchers who has given his considerations to ponder asymptotic behavior of solutions to the equations proposed to study different flexible materials with voids (Cowin, Goodman and Nunziato \cite{1,2,3,4}), which have decent physical properties, are utilized broadly in engineering, such as vehicles, aero planes, expansive space structures and so on. Due to their broad applications,  some of this interest established by many researchers comes from the need to establish results concerning existence and stabilization the elasticity problems.\newline
In addition to the conventional elastic effects, materials with voids have a microstructure in which the mass at each place is calculated by multiplying the material matrix  mass density by the volume fraction. Nunziato and Cowin \cite{4} pioneered the latter concept in their groundbreaking work on elastic materials with voids. Iesan \cite{4.1,4.2, 4.32} and Iesan an Quintanilla \cite{4.42} expanded the hypothesis by including temperature and microtemperatures.\newline
According to our knowledge, the evaluation of the temporal decay in one-dimensional porouselastic
substances become pioneered with the aid  Quintanilla \cite{4.81} where he proved that porous-viscosity become not
robust sufficient to exponentially stabilize the system. Interestingly, Casas and Quintanilla \cite{4.82}
proved that the mixture of porous-viscosity and temperature additionally lacks exponential stability.
However, the identical authors \cite{4.83} confirmed that the mixture of porous-viscosity and thermal
effects (each temperature and microtemparatures) stabilized the system exponentially. Similarly,
Magana and Quintanilla \cite{4.84} proved that viscoelasticity collectively with microtemperatures produced
exponential stability, while viscoelasticity collectively with temperature lacks exponential
stability. \newline 
It is natural to think that porous-elastic system with dissipation due to only microtemperatures will definitely lack exponential stability. However T. Apalara \cite{7} establish the contrary, he proved that a porous system with microtempearture decays exponentially if and only if $\chi=0$ where $\chi=\frac{\mu}{\rho}-\frac{\delta}{J}$, otherwise the system is polynomially stable. \newline 
The equations for one dimensional  porous elastic system with microtempeartureare of the form
\begin{equation}\label{0.1}
\left\{
\begin{aligned}
&\r u_{tt}=T_{x},\\
 &J\phi_{tt}=H_{x}+G,\\
 &\r E_{t}=P_{x}+q-Q.
 \end{aligned}
 \right.
 \end{equation}
   Where  $(x,t)\in (0,l)\times (0,\infty)$, $t$ is the time, $x$ is the distance along the center line of 
   the beam structure and $l$ is the length of the beam, $T$ is the stress, $H$ is the equilibrated
   stress, $G$ is the equilibrated body force, $q$ is the heat flux vector, $P$ is the first heat flux moment, $Q$ is the mean heat flux and $E$ is the first moment of energy. The functions $u(x,t)$ and $\phi(x,t)$ are  respectively, the displacement of the solid elastic material and the volume fraction.
   The constitutive equations are given by 
   \begin{equation}\label{0.2}
   \left\{
   \begin{aligned}
            &  T=\mu u_{x}+b\phi,\\
              &  H=\delta \phi_{x}-dw,\\                                                     
               &  G=-bu_{x}-\xi\phi,\\
               & \r E=-\a w-d\phi_{x},\\
               & P=-\kappa w_{x},\\
               &q=k_{1}w,\\
               & Q=k_{2}w.
   \end{aligned}
   \right.
   \end{equation}     
   where $w$ is the microtemperature, $k_{1}$, $k_{2}$, $d$,  $\mu$, $\delta$, $\a$, $\kappa$ and $\xi$ are positive constants such that  $\mu\xi-b^{2}>0$.      
Substitute system $\eqref{0.2}$ in $\eqref{0.1}$ we get the porous elastic system with microtemperature 
\begin{equation}\label{0.3}
\r u_{tt}-\mu u_{xx}- b \phi_{x}=0, \quad \mbox{in} \quad (0,l)\times (0.\infty),
\end{equation}
\begin{equation}\label{0.4}
J\phi_{tt}-\delta\phi_{xx}+b u_{x}+\xi\phi+dw_{x}=0, \quad \mbox{in} \quad (0,l)\times (0.\infty),
\end{equation}
\begin{equation}\label{0.5}
\a w_{t}-\kappa w_{xx}+d\phi_{tx}+k w=0, \quad \mbox{in} \quad (0,l)\times (0.\infty)
\end{equation}
where $k=k_{1}-k_{2}>0$. \newline 
In this paper we consider the following porous system with second spectrum free and microtemperature: 
\begin{equation}\label{1.1}
\r u_{tt}-\mu u_{xx}- b \phi_{x}=0, \quad \mbox{in} \quad (0,l)\times (0.\infty),
\end{equation}
\begin{equation}\label{1.2}
-Ju_{ttx}-\delta \phi_{xx}+b u_{x}+\xi\phi+d w_{x}=0, \quad \mbox{in} \quad (0,l)\times (0.\infty),
\end{equation}
\begin{equation}\label{1.3}
\a w_{t}-\kappa w_{xx}+d\phi_{tx}+k w=0, \quad \mbox{in} \quad (0,l)\times (0.\infty),
\end{equation}
With Dirichlet boundary conditions 
\begin{equation}\label{1.4}
u(0,t)=u(l,t)=\phi(0,t)=\phi(l,t)=w(0,t)=w(l,t)=0, \quad t>0,
\end{equation}
and the initial conditions are 
\begin{equation}\label{1.5}
u(x,0)=u_{0}(x), \quad u_{t}(x,0)=u_{1}(x), \quad u_{tt}(x,0)=u_{2}(x), \quad \phi(x,0)=\phi_{0}(x),\quad w(x,0)=w_{0}(x) \quad x\in (0,l).
\end{equation}
This system \eqref{1.1}-\eqref{1.5} is obtained by following the procedure of Elishakoff \cite{26} which involves replacing the term $\phi_{tt}$  in \ref{0.4}  by $-u_{xtt}$  based on d' Alembert's  principle for dynamic equilibrium. This eliminates the second
spectrum of frequency and its damaging consequences for wave propagation speed. The goal of this work is to get exponential decay without assuming any conditions on the physical parameters.
The dissipation of system $\eqref{1.1}-\eqref{1.5}$ is obtained from the definition of energy. Indeed, multiply equation $\eqref{1.1}$ by $u_{t}$, integrate by parts over $(0,l)$ and using boundary conditions $\eqref{1.4}$ we get 
\begin{equation}\label{0.11}
\frac{1}{2}\frac{d}{dt}\left(\r||u_{t}||^{2}+\mu||u_{x}||^{2}\right)+b\int_{0}^{l}\phi u_{xt}dx=0.
\end{equation}
Multiply equation $\eqref{1.2}$ by $\phi_{t}$, integrate by parts over $(0,l)$ and using boundary conditions $\eqref{1.4}$ we get 
\begin{equation}\label{0.12}
\frac{1}{2}\frac{d}{dt}\left(\delta||\phi_{x}||^{2}+\xi||\phi||^{2}\right)+b\int_{0}^{l}\phi_{t}u_{x}dx+J\int_{0}^{l}u_{tt}\phi_{tx}dx+d\int_{0}^{l}w_{x}\phi_{t}dx=0.
\end{equation}
From equation $\eqref{1.1}$ we obtain that $$\phi_{xt}=\frac{\r}{b}u_{ttt}-\frac{\mu}{b}u_{xxt},$$
substitute $\phi_{xt}$ in equation $\eqref{0.12}$ we arrive at 
\begin{equation}\label{0.13}
\frac{1}{2}\frac{d}{dt}\left(\delta||\phi_{x}||^{2}+\xi||\phi||^{2}+\frac{J\r}{b}\|u_{tt}\|^{2}+\frac{J\mu}{b}\|u_{xt}\|^{2}\right)+b\int_{0}^{l}\phi_{t}u_{x}dx+d\int_{0}^{l}w_{x}\phi_{t}dx=0.
\end{equation}
Now multiply equation $\eqref{1.3}$ by $w$, integrate by parts over $(0,l)$ and using boundary conditions $\eqref{1.4}$ we get 
\begin{equation}\label{0.14}
\frac{1}{2}\frac{d}{dt}\left(\a||w||^{2}\right)+\kappa\int_{0}^{l}w_{x}^{2}+d\int_{0}^{l}\phi_{tx}w+k\int_{0}^{l}w^{2}=0.
\end{equation}
Add equations $\eqref{0.11}$, $\eqref{0.13}$ and $\eqref{0.14}$ we get
\begin{equation*}
\frac{1}{2}\frac{d}{dt}\left(\r||u_{t}||^{2}+\mu||u_{x}||^{2}+\delta||\phi_{x}||^{2}+\xi||\phi||^{2}+\frac{J\r}{b}\|u_{tt}\|^{2}+\frac{J\mu}{b}\|u_{xt}\|^{2}+\a||w||^{2}+2b(u_{x},\phi) \right)=-\kappa\int_{0}^{l}w_{x}^{2}-k\int_{0}^{l}w^{2}.
\end{equation*}
Define 
$$E(t)=\frac{1}{2}\left(\r||u_{t}||^{2}+\mu||u_{x}||^{2}+\delta||\phi_{x}||^{2}+\xi||\phi||^{2}+\frac{J\r}{b}\|u_{tt}\|^{2}+\frac{J\mu}{b}\|u_{xt}\|^{2}+\a||w||^{2}+2b(u_{x},\phi) \right).$$
Add then subtract $\frac{b^{2}}{2\xi}\|u_{x}\|^{2}$ to the right side of the above equation we arrive at
$$E(t)=\frac{1}{2}\left(\r||u_{t}||^{2}+\frac{J\r}{b}||u_{tt}||^{2}+(\mu-b^{2}/\xi)||u_{x}||^{2}+\frac{J\mu}{b}||u_{xt}||^{2}+\delta||\phi_{x}||^{2}+||\frac{b}{\sqrt{\xi}}u_{x}+\sqrt{\xi}\phi||^{2}+\a||w||^{2}\right)$$
Where $||.||$ denotes the $L^{2}$--norm. \newline 
$E(t)$ preserves its positivity for $\mu\xi-b^{2}>0$
and the dissipation law is given by 
$$ \frac{d}{dt}E(t)=-\kappa\int_{0}^{l}w_{x}^{2}-k\int_{0}^{l}w^{2}.$$

\section{Well-Posedness}
\SE
The aim of this section is to show the existence and uniqueness of weak solution of system $\eqref{1.1}-\eqref{1.5}$. For that reason we will use the classical Faedo-Galerkin approximation along with an a priori estimates then passing through the limits using compactness arguments.\newline
Define the Hilbert space 
$$\mathcal{H}:=H_{0}^{1}(0,l)\times H_{0}^{1}(0,l)\times L^{2}(0,l)\times H_{0}^{1}(0,l)\times L^{2}(0,l) .$$ 
Now multiply the equations $\eqref{1.1}$, $\eqref{1.2}$ and $\eqref{1.3}$ by $\bar{u}$, $\bar{\phi}$, $\bar{w} \in H_{0}^{1}(0,l)$ respectively and integrate by parts over $(0,l)$ we get using boundary conditions $\eqref{1.4}$ 
\begin{equation}\label{W2.1}
\left\{
\begin{aligned}
&\r(u_{tt},\bar{u})+\mu(u_{x},\bar{u}_{x})+b(\phi,\bar{u}_{x})=0\\
&J(u_{tt},\bar{\phi}_{x})+\delta(\phi_{x},\bar{\phi}_{x})+b(u_{x},\bar{\phi})+\xi(\phi, \bar{\phi})+d(w_{x},\bar{\phi}) =0\\
&\a (w_{t},\bar{w})+\kappa( w_{x},\bar{w}_{x})+d(\phi_{tx},\bar{w}) +k( w, \bar{w})=0.
\end{aligned}
\right.
\end{equation}
\begin{definition}
Let the initial data $(u_{0},u_{1},u_{2},\phi_{0}, w_{0})\in \mathcal{H}$ then a function $V=(u,u_{t},u_{tt},\phi, w)\in C(0,T;\mathcal{H})$ is said to be a weak solution of $\eqref{1.1}-\eqref{1.5}$ if it is a solution of the weak problem $\eqref{W2.1}$ for almost $t\in [0,T].$
\end{definition}
\begin{theorem}
Suppose that the initial data $(u_{0},u_{1},u_{2},\phi_{0}, w_{0})\in \mathcal{H}$ then system $\eqref{1.1}-\eqref{1.5}$ have a weak solution satisfying 
$$u \in  L^{\infty}\left(0,T;H_{0}^{1}(0,l)\right), \quad  u_{t} \in L^{\infty}\left(0,T;H_{0}^{1}(0,l)\right),$$
$$u_{tt} \in  L^{\infty}\left(0,T;L^{2}(0,l)\right), \quad \phi \in  L^{\infty}\left(0,T;H_{0}^{1}(0,l)\right),$$
$$w \in L^{\infty}\left(0,T;L^{2}(0,l)\right)\cap L^{2}\left(0,T;H_{0}^{1}(0,l)\right),$$
where the solution $V=(u,u_{t},u_{tt},\phi, w)$ depends continuously on the initial data in $\mathcal{H}$. In particular $V$ is unique solution of system $\eqref{1.1}-\eqref{1.5}$.
\end{theorem}
\noindent{\textbf{Proof.}}
We will use the Faedo-Galerkin method to prove the above theorem and we proceed in five steps.
\textbf{Step 1. Approximated solution}
 Let $(u_{0},u_{1},u_{2},\phi_{0}, w_{0})\in \mathcal{H}$.
Let $\{\eta_{i}\}_{i=1}^{\infty} \subset C^{\infty}([0,l])$ be basis for $H_{0}^{1}(0,l)$, and let $V^{m}=span\{\eta_{i}\}_{i=1}^{m}$. Now we introduce 
\begin{equation}\label{W2.2}
u^{m}=\sum_{i=0}^{m}a_{i}(t)\eta_{i}(x), \quad \phi^{m}=\sum_{i=0}^{m}b_{i}(t)\eta_{i}(x),\quad w^{m}=\sum_{i=0}^{m}c_{i}(t)\eta_{i}(x),
\end{equation}
which solves the following approximated problem for $\bar{u}$, $\bar{\phi}$, $\bar{w} \in V^{m}$ 
\begin{equation}\label{W2.3}
\left\{
\begin{aligned}
&\r(u_{tt}^{m},\bar{u})+\mu(u_{x}^{m},\bar{u}_{x})+b(\phi^{m},\bar{u}_{x})=0\\
&J(u_{tt}^{m},\bar{\phi}_{x})+\delta(\phi_{x}^{m},\bar{\phi}_{x})+b(u_{x}^{m},\bar{\phi})+\xi(\phi^{m}, \bar{\phi})+d(w_{x}^{m},\bar{\phi}) =0\\
&\a (w_{t}^{m},\bar{w})+\kappa( w_{x}^{m},\bar{w}_{x})+d(\phi_{tx}^{m},\bar{w}) +k( w^{m}, \bar{w})=0.
\end{aligned}
\right.
\end{equation}
with initial conditions $$\left(u^{m}(0),u_{t}^{m}(0),u_{tt}^{m}(0),\phi^{m}(0), w^{m}(0)\right)=(u_{0}^{m},u_{1}^{m},u_{2}^{m},\phi_{0}^{m}, w_{0}^{m})$$
such that 
 $$(u_{0}^{m},u_{1}^{m},u_{2}^{m},\phi_{0}^{m}, w_{0}^{m})\to (u_{0},u_{1},u_{2},\phi_{0}, w_{0}) \quad \mbox{strongly in} \quad \mathcal {H}.$$ 
 By using the Carathoedory theorem for standard ordinary differential equations theory, system $\eqref{W2.3}$  has a local solution $\left(u^{m}(t),u_{t}^{m}(t),u_{tt}^{m}(t),\phi^{m}(t), w^{m}(t)\right)$  on the maximal interval $[0,t_{m})$ with $0<t_{m}\leq T$ for every $m\in \mathbb{N}$ \newline 
\noindent\textbf{Step 2. A priori estimates.}\newline 
 Let $\bar{u}=u_{t}^{m}$, $\bar{\phi}=\phi_{t}^{m}$ and $\bar{w}=w^{m}$ and taking into consideration from equation $\eqref{1.1}$ that $$\phi_{xt}=\frac{\r}{b}u_{ttt}-\frac{\mu}{b}u_{xxt},$$ then $$\phi_{xt}^{m}=\frac{\r}{b}u_{ttt}^{m}-\frac{\mu}{b}u_{xxt}^{m},$$ system $\eqref{W2.3}$ becomes 
\begin{equation}\label{W2.4}
\left\{
\begin{aligned}
&\r(u_{tt}^{m},u_{t}^{m})+\mu(u_{x}^{m},u_{tx}^{m})+b(\phi^{m},u_{tx}^{m})=0\\
&\frac{J\r}{b}\left(u_{tt}^{m},u_{ttt}^{m}\right)+\frac{J\mu}{b}\left(u_{ttx}^{m},u_{tx}^{m}\right)+\delta(\phi_{x}^{m},\phi_{tx}^{m})+b(u_{x}^{m},\phi_{t}^{m})+\xi(\phi^{m}, \phi_{t}^{m})+d(w_{x}^{m},\phi_{t}^{m}) =0\\
&\a (w_{t}^{m},w^{m})+\kappa( w_{x}^{m},w_{x}^{m})+d(\phi_{tx}^{m},w^{m}) +k( w^{m}, w^{m})=0.
\end{aligned}
\right.
\end{equation}
Which is equivalent to 
\begin{equation}\label{W2.5}
\left\{
\begin{aligned}
&\frac{d}{2dt}\left(\r||u_{t}^{m}||^{2}\right)+\frac{d}{2dt}\left(\mu||u_{x}^{m}||^{2}\right)+b(\phi^{m},u_{tx}^{m})=0\\
&\frac{d}{2dt}\left(\frac{J\r}{b}\|u_{tt}^{m}\|^{2}\right)+\frac{d}{2dt}\left(\frac{J\mu}{b}\|u_{tx}^{m}\|^{2}\right)+\frac{d}{2dt}\left(\delta||\phi_{x}^{m}||^{2}\right)+b(u_{x}^{m},\phi_{t}^{m})+\frac{d}{2dt}\left(\xi||\phi^{m}||^{2}\right)+d(w_{x}^{m},\phi_{t}^{m})=0\\
&\frac{d}{2dt}\left(\a||w^{m}||^{2}\right)+\kappa||w_{x}^{m}||^{2}+d(\phi_{tx}^{m},w^{m}) +k||w^{m}||^{2}=0,
\end{aligned}
\right.
\end{equation}
where $||.||$ denotes the norm in $L^{2}(0,l)$.\newline
Add the above two equations we get 
\begin{equation}\label{W2.6}
\begin{aligned}
&\frac{d}{2dt}\left(\r||u_{t}^{m}||^{2}+\mu||u_{x}^{m}||^{2}+\frac{J\r}{b}\|u_{tt}^{m}\|^{2}+\frac{J\mu}{b}\|u_{tx}^{m}\|^{2}+\delta||\phi_{x}^{m}||^{2}+\xi||\phi^{m}||^{2}+2b(u_{x}^{m},\phi^{m})+\a||w^{m}||^{2}\right)\\
&+\kappa||w_{x}^{m}||^{2}+k||w^{m}||^{2} =0.
\end{aligned}
\end{equation}
Let $$E^{m}(t)=\frac{1}{2}\left(\r||u_{t}^{m}||^{2}+\mu||u_{x}^{m}||^{2}+\frac{J\r}{b}\|u_{tt}^{m}\|^{2}+\frac{J\mu}{b}\|u_{tx}^{m}\|^{2}+\delta||\phi_{x}^{m}||^{2}+\xi||\phi^{m}||^{2}+2b(u_{x}^{m},\phi^{m})+\a||w^{m}||^{2}\right).$$
Then equation $\eqref{W2.6}$ becomes 
\begin{equation}\label{W2.7}
\frac{d}{dt}E^{m}(t)+\kappa||w_{x}^{m}||^{2}+k||w^{m}||^{2} =0. 
\end{equation}
Now integrate $\eqref{W2.7}$ from $0$ to $t<t_{m}$, we obtain from the choice of the initial data that for all $t\in [0,T]$  and for every $m\in \mathbb{N}$  that 
\begin{equation}\label{W2.8}
E^{m}(t)+\kappa \int_{0}^{t}||w_{x}^{m}(s)||^{2}ds+k\int_{0}^{t}||w^{m}(s)||^{2}ds\leq C_{0},
\end{equation}
where $C_{0}$ is a positive constant depending on the initial data.\newline 

 \noindent\textbf{Step 3. Passing to the limit.}
Using $\eqref{W2.8}$  and by the definition of $E^{m}(t)$  we obtain that 
\begin{equation*}
\left\{
\begin{aligned}
&\{u^{m}\} \quad \mbox{is bounded in} \quad L^{\infty}\left(0,T;H_{0}^{1}(0,l)\right)\\
&\{u_{t}^{m}\} \quad \mbox{is bounded in} \quad L^{\infty}\left(0,T;H_{0}^{1}(0,l)\right)\\
&\{u_{tt}^{m}\} \quad \mbox{is bounded in} \quad L^{\infty}\left(0,T;L^{2}(0,l)\right)\\
&\{\phi^{m}\} \quad \mbox{is bounded in} \quad L^{\infty}\left(0,T;H_{0}^{1}(0,l)\right)\\
&\{w^{m}\} \quad \mbox{is bounded in} \quad  L^{\infty}\left(0,T;L^{2}(0,l)\right)\cap L^{2}\left(0,T;H_{0}^{1}(0,l)\right)
\end{aligned}
\right.
\end{equation*}
Then we can extract a subsequence of $\{u^{m}\}$, $\{\phi^{m}\}$ and $\{w^{m}\}$ and still denoted by $\{u^{m}\}$, $\{\phi^{m}\}$ and $\{w^{m}\}$, such that  
\begin{equation*}
\left\{
\begin{aligned}
&u^{m}\to u \quad \mbox{weakly star in} \quad L^{\infty}\left(0,T;H_{0}^{1}(0,l)\right)\\
&u_{t}^{m}\to u_{t} \quad \mbox{weakly star in} \quad L^{\infty}\left(0,T;H_{0}^{1}(0,l)\right)\\
&u_{tt}^{m}\to u_{tt} \quad \mbox{weakly star in} \quad L^{\infty}\left(0,T;L^{2}(0,l)\right)\\
&\phi^{m}\to \phi \quad \mbox{weakly star in} \quad L^{\infty}\left(0,T;H_{0}^{1}(0,l)\right)\\
&w^{m}\to w \quad \mbox{weakly star in} \quad L^{\infty}\left(0,T;L^{2}(0,l)\right)\\
&w^{m}\to w \quad \mbox{weakly in} \quad L^{2}\left(0,T;H_{0}^{1}(0,l)\right)
\end{aligned}
\right.
\end{equation*}
Now pass to the limits in the approximate variational problem $\eqref{W2.3}$ we get a weak solution satisfying 
$$u \in  L^{\infty}\left(0,T;H_{0}^{1}(0,l)\right), \quad  u_{t} \in L^{\infty}\left(0,T;H_{0}^{1}(0,l)\right),$$
$$u_{tt} \in  L^{\infty}\left(0,T;L^{2}(0,l)\right), \quad \phi \in  L^{\infty}\left(0,T;H_{0}^{1}(0,l)\right),$$
$$w \in L^{\infty}\left(0,T;L^{2}(0,l)\right)\cap L^{2}\left(0,T;H_{0}^{1}(0,l)\right).$$
\noindent\textbf{Step 4. Initial data.}
Knowing that $$H_{0}^{1}(0,l)\subset L^{2}(0,l)\subset H^{-1}(0,l),$$ 
where $H^{-1}(0,l)$ is the dual space of $H_{0}^{1}(0,l)$. \newline 
By using Aubin-Lions lemma, see \cite{47}, we obtain that 
$ L^{\infty}\left(0,T;H_{0}^{1}(0,l)\right)$ is compactly embedded in $C\left(0,T;L^{2}(0,l)\right)$. This implies that 
$$u^{m}\to u \quad \mbox{strongly in} \quad C\left(0,T;L^{2}(0,l)\right),$$
$$u_{t}^{m}\to u_{t} \quad \mbox{strongly in} \quad C\left(0,T;L^{2}(0,l)\right),$$

Hence, $$\left(u(0),u_{t}(0)\right)=(u_{0},u_{1}).$$
Now differentiate with respect to $t$ the first equation of system $\eqref{W2.3}$ we get 
$$\r(u_{ttt}^{m},\bar{u})+\mu(u_{tx}^{m},\bar{u}_{x})+b(\phi_{t}^{m},\bar{u}_{x})=0,$$ 
for all $\bar{u}\in H_{0}^{1}(0,l)$.\newline 
Multiply the above equation  by a test function 
$$\lambda\in H_{0}^{1}(0,T), \quad \mbox{such that} \quad \lambda(0)=1, \quad \lambda(T)=0,$$
and then integrate by parts over $[0,T]$ 
$$-\r(u_{2}^{m},\bar{u})-\r\int_{0}^{T}(u_{tt}^{m},\bar{u})\lambda_{t}dt+\mu\int_{0}^{T}(u_{tx}^{m},\bar{u}_{x})\lambda dt+b\int_{0}^{T}(\phi_{t}^{m},\bar{u}_{x})\lambda dt=0.$$
Take the limit $m\to\infty$, we arrive at 
\begin{equation}\label{18}
-\r(u_{2},\bar{u})-\r\int_{0}^{T}(u_{tt},\bar{u})\lambda_{t}dt+\mu\int_{0}^{T}(u_{tx},\bar{u}_{x})\lambda dt+b\int_{0}^{T}(\phi_{t},\bar{u}_{x})\lambda dt=0.
\end{equation}
Now differentiate the first equation of system $\eqref{W2.1}$ with respect to time, then multiply the result by $\lambda$ under the same conditions above and integrate by parts over $[0,T]$ we get 
\begin{equation}\label{19}
-\r(u_{tt}(0),\bar{u})-\r\int_{0}^{T}(u_{tt},\bar{u})\lambda_{t}dt+\mu\int_{0}^{T}(u_{tx},\bar{u}_{x})\lambda dt+b\int_{0}^{T}(\phi_{t},\bar{u}_{x})\lambda dt=0.
\end{equation}
Combine the two equations $\eqref{18}$ and $\eqref{19}$ we obtain that $u_{tt}(0)=u_{2}.$ In the same way we can get that $(\phi(0), w(0))=(\phi_{0}, w_{0})$.\newline
\noindent\textbf{Step 5. Continuous dependence on initial data.} Let $V_{1}(t)=(u,u_{t},u_{tt},\phi, w)$ and $V_{2}(t)=(\tilde{u},\tilde{u}_{t},\tilde{u}_{tt},\tilde{\phi}, \tilde{w})$ be two solutions of the system $\eqref{1.1}-\eqref{1.3}$ with initial data $V_{1}(0)=(u_{0},u_{1},u_{2},\phi_{0}, w_{0})$ and $V_{2}(0)=(\tilde{u}_{0},\tilde{u}_{1},\tilde{u}_{1},\tilde{\phi}_{0}, \tilde{w}_{0})$ such that $V_{1}(0),V_{2}(0)\in \mathcal{H}$. Then $(U,U_{t},U_{tt},\Phi, W)=V_{1}(t)-V_{2}(t)$ satisfies the following equations 
\begin{equation}\label{U1.1}
\r U_{tt}-\mu U_{xx}- b \Phi_{x}=0, \quad \mbox{in} \quad (0,l)\times (0.\infty),
\end{equation}
\begin{equation}\label{U1.2}
-JU_{ttx}-\delta\Phi_{xx}+b U_{x}+\xi\Phi+d W_{x}=0, \quad \mbox{in} \quad (0,l)\times (0.\infty),
\end{equation}
\begin{equation}\label{U1.3}
\a W_{t}-\kappa W_{xx}+d\Phi_{tx}+k W=0, \quad \mbox{in} \quad (0,l)\times (0.\infty),
\end{equation}

with initial data $(U_{0},U_{1},U_{2},\Phi_{0}, W_{0})=V_{1}(0)-V_{2}(0)$.\newline 
Now multiply $\eqref{U1.1}$ by $U_{t}$, $\eqref{U1.2}$ by $\Phi_{t}$ and $\eqref{U1.3}$ by $W$ then integrate the result over $(0,l)$ we arrive at 
\begin{equation}\label{27}
\frac{d}{dt}\hat{E}(t)=-\kappa\int_{0}^{l}W_{x}^{2}dx-k\int_{0}^{l}W^{2}dx,
\end{equation}
where $\hat{E}(t)$ is the energy related to $V_{1}(t)-V_{2}(t)$ and defined by 
$$\hat{E}(t)=\frac{1}{2}\left(\r||U_{t}||^{2}+\frac{J\r}{b}||U_{tt}||^{2}+(\mu-b^{2}/\xi)||U_{x}||^{2}+\frac{J\mu}{b}||U_{xt}||^{2}+\delta||\Phi_{x}||^{2}+||\frac{b}{\sqrt{\xi}}U_{x}+\sqrt{\xi}\Phi||^{2}+\a||W||^{2}\right).$$
Integrate $\eqref{27}$ over $(0,t)$, w get that there exists a positive constant $C_{T}$ such that for any $t\in [0,T]$,
$$\hat{E}(t)\leq C_{T}\hat{E}(0),$$
which implies that the weak solution depend continuously on the initial data. Consequently the weak solution of system $\eqref{1.1}-\eqref{1.3}$ is unique.

\section{Exponential stability}
\SE
\begin{theorem}\label{Thm1}
The energy $E(t)$ of the system $\eqref{1.1}-\eqref{1.5}$ decays exponentially as time t tends to infinity. That is, there
exist two positive constants $M$ and $\omega$  independent of the initial data and independent of any relationship between
coefficients such that
$$E(t)\leq ME(0)e^{-\omega t}, \quad \forall t\geq 0.$$
\end{theorem}
The proof of \textbf{Theorem 3} will be established through two lemmas. First, we set
$$\mathcal{F}(t)=\r\int_{0}^{l}u_{t}udx+\frac{J\mu}{b}\int_{0}^{l}u_{xt}u_{x}dx.$$
\begin{lemma}\label{lemma 2}
Let $(u,\phi, w)$ be a solution of the system $\eqref{1.1}-\eqref{1.5}$. Then we have 
\begin{equation}\label{2.1}
\begin{aligned}
\frac{d}{dt}\mathcal{F}(t)&\leq\left(\frac{J\mu}{b}+2\r c_{p}\right)\int_{0}^{l}|u_{xt}|^{2}dx -\r\int_{0}^{l}|u_{t}|^{2}dx-(\mu-b^{2}/\xi)\int_{0}^{l}|u_{x}|^{2}dx-\frac{\delta}{2}\int_{0}^{L}|\phi_{x}|^{2}dx-\frac{J\r}{b}\int_{0}^{l}|u_{tt}|^{2}dx\\
&-\int_{0}^{l}\left|\frac{b}{\sqrt{\xi}}u_{x}+\sqrt{\xi}\phi\right|^{2}dx+\frac{d^{2}c_{p}}{2\delta}\int_{0}^{l}|w_{x}|^{2}dx.
\end{aligned}
\end{equation}
where $c_{p}$ is the Poincare's constant.
\end{lemma}
\noindent\textbf{Proof:}
Multiply equation $\eqref{1.1}$ by $u$ and integrate by parts over $(0,l)$ we get 
$$\r\int_{0}^{l} u_{tt}udx+\mu\int_{0}^{l} |u_{x}|^{2}dx + b \int_{0}^{l}u_{x}\phi dx=0,$$
add then subtract the term $\frac{b}{\sqrt{\xi}}\int_{0}^{l}|u_{x}|^{2}$ from the above equation we obtain 
$$\r\int_{0}^{l} u_{tt}udx+(\mu-b^{2}/\xi)\int_{0}^{l} |u_{x}|^{2} dx+\frac{b}{\sqrt{\xi}} \int_{0}^{l}\left(\frac{b}{\sqrt{\xi}}u_{x}+\sqrt{\xi}\phi\right)u_{x}dx=0,$$
taking into account that $\frac{d}{dt}(u_{t}u)=u_{tt}u+|u_{t}|^{2}$ we arrive at 
\begin{equation}\label{2.2}
\frac{d}{dt}\left ( \r\int_{0}^{l}u_{t}u dx\right )=\r\int_{0}^{l}|u_{t}|^{2}dx-(\mu-b^{2}/\xi)\int_{0}^{l} |u_{x}|^{2}dx -\frac{b}{\sqrt{\xi}} \int_{0}^{l}\left(\frac{b}{\sqrt{\xi}}u_{x}+\sqrt{\xi}\phi\right)u_{x}dx.
\end{equation}
Multiply equation $\eqref{1.2}$ by $\phi$ and integrate by parts over $(0,l)$ we get 
\begin{equation}\label{2.3}
J\int_{0}^{l}u_{tt}\phi_{x}dx+\delta \int_{0}^{l}|\phi_{x}|^{2}dx+\sqrt{\xi} \int_{0}^{l}\left(\frac{b}{\sqrt{\xi}}u_{x}+\sqrt{\xi}\phi\right)\phi dx+d\int_{0}^{l}w_{x}\phi dx=0
\end{equation}
From equation $\eqref{1.1}$ we get that $\phi_{x}=\frac{\r}{b}u_{tt}-\frac{\mu}{b}u_{xx},$ then substitute $\phi_{x}$ in equation $\eqref{2.3}$ and taking into account that $\frac{d}{dt}(u_{tx}u_{x})=u_{ttx}u_{x}+|u_{tx}|^{2}$ we obtain
$$\frac{d}{dt}\left(\frac{J\mu}{b}\int_{0}^{l}u_{xt}u_{x}dx\right)-\frac{J\mu}{b}\int_{0}^{l}|u_{tx}|^{2}dx+\frac{J\r}{b}\int_{0}^{l}|u_{tt}|^{2}dx+\delta \int_{0}^{l}|\phi_{x}|^{2}dx+\sqrt{\xi} \int_{0}^{l}\left(\frac{b}{\sqrt{\xi}}u_{x}+\sqrt{\xi}\phi\right)\phi dx+d\int_{0}^{l}w_{x}\phi dx=0,$$
Using Poincare's and Young's inequality we get 
\begin{equation}\label{2.35}
\begin{aligned}
\frac{d}{dt}\left(\frac{J\mu}{b}\int_{0}^{l}u_{xt}u_{x}dx\right)&\leq\frac{J\mu}{b}\int_{0}^{l}|u_{tx}|^{2}dx-\frac{J\r}{b}\int_{0}^{l}|u_{tt}|^{2}dx-\frac{\delta}{2} \int_{0}^{l}|\phi_{x}|^{2}dx\\
&-\sqrt{\xi} \int_{0}^{l}\left(\frac{b}{\sqrt{\xi}}u_{x}+\sqrt{\xi}\phi\right)\phi dx+\frac{d^{2}c_{p}}{2\delta}\int_{0}^{l}|w_{x}|^{2}dx.
\end{aligned}
\end{equation}
Add the two equations $\eqref{2.2}$ and $\eqref{2.35}$ we obtain 
\begin{equation*}
\begin{aligned}
\frac{d}{dt}\left ( \r\int_{0}^{l}u_{t}u dx+\frac{J\mu}{b}\int_{0}^{l}u_{xt}u_{x}dx\right) &\leq\frac{J\mu}{b}\int_{0}^{l}|u_{xt}|^{2}dx +\r\int_{0}^{l}|u_{t}|^{2}dx-(\mu-b^{2}/\xi)\int_{0}^{l}|u_{x}|^{2}dx-\frac{\delta}{2}\int_{0}^{L}|\phi_{x}|^{2}dx\\
&-\frac{J\r}{b}\int_{0}^{l}|u_{tt}|^{2}dx-\int_{0}^{l}\left|\frac{b}{\sqrt{\xi}}u_{x}+\sqrt{\xi}\phi\right|^{2}dx+\frac{d^{2}c_{p}}{2\delta}\int_{0}^{l}|w_{x}|^{2}dx.
\end{aligned}
\end{equation*}
Add and subtract the term $\r\int_{0}^{l}|u_{t}|^{2}dx$ to the right side of the above inequality then use Poincare's inequality we get the desired result. $\hfill{\blacksquare}$ \newline
Set $$\mathcal{G}(t)=-J\int_{0}^{l}u_{tx}\left(\frac{b}{\sqrt{\xi}}u_{x}+\sqrt{\xi}\phi\right)dx-\delta\frac{\r b}{\mu \sqrt{\xi}}\int_{0}^{l}u_{tx}\phi dx+\frac{\a C_{1}}{d}\int_{0}^{l}w u_{t}dx.$$
\begin{lemma}\label{lemma 3}
Let $(u,\phi, w)$ be a solution of the system $\eqref{1.1}-\eqref{1.5}$. Then we have 
\begin{equation}\label{2.1}
\begin{aligned}
\frac{d}{dt}\mathcal{G}(t)&\leq-\frac{Jb}{2\sqrt{\xi}}\int_{0}^{l}|u_{xt}|^{2}dx +\frac{\eps_{3}}{2}\int_{0}^{l}|u_{t}|^{2}dx+\frac{\eps_{1}}{2}\int_{0}^{l}|u_{tt}|^{2}dx-\delta(\mu-b^{2}/\xi)\frac{\sqrt{\xi}}{\mu}\int_{0}^{L}|\phi_{x}|^{2}dx\\
&-\frac{\sqrt{\xi}}{2}\int_{0}^{l}\left|\frac{b}{\sqrt{\xi}}u_{x}+\sqrt{\xi}\phi\right|^{2}dx+C_{2}\int_{0}^{l}|w|^{2}dx+C_{3}\int_{0}^{l}|w_{x}|^{2}dx.
\end{aligned}
\end{equation}
where $C_{1}$, $C_{2}$ and $C_{3}$ are positive constant to be determined.
\end{lemma}
\noindent{\textbf{Proof.}}
Multiply equation $\eqref{1.2}$ by $\left(\frac{b}{\sqrt{\xi}}u_{x}+\sqrt{\xi}\phi\right)$ we get:
\begin{equation}\label{2.7}
\begin{aligned}
&-J\int_{0}^{l}u_{ttx}\left(\frac{b}{\sqrt{\xi}}u_{x}+\sqrt{\xi}\phi\right)dx+\delta\int_{0}^{l}\phi_{x}\left(\frac{b}{\sqrt{\xi}}u_{x} +\sqrt{\xi}\phi\right)_{x}dx\\
&=-\sqrt{\xi}\int_{0}^{l}\left|\frac{b}{\sqrt{\xi}}u_{x}+\sqrt{\xi}\phi\right|^{2}dx-d\int_{0}^{l}w_{x}\left(\frac{b}{\sqrt{\xi}}u_{x}+\sqrt{\xi}\phi\right)dx.
\end{aligned}
\end{equation}
Using Young's inequality we obtain 
\begin{equation}\label{2.8}
\begin{aligned}
&-J\int_{0}^{l}u_{ttx}\left(\frac{b}{\sqrt{\xi}}u_{x}+\sqrt{\xi}\phi\right)dx+\delta\int_{0}^{l}\phi_{x}\left(\frac{b}{\sqrt{\xi}}u_{x} +\sqrt{\xi}\phi\right)_{x}dx\\
&=-\frac{\sqrt{\xi}}{2}\int_{0}^{l}\left|\frac{b}{\sqrt{\xi}}u_{x}+\sqrt{\xi}\phi\right|^{2}dx+\frac{d^{2}}{2\sqrt{\xi}}\int_{0}^{l}|w_{x}|^{2}dx.
\end{aligned}
\end{equation}
Add then subtract the term $\frac{\mu\xi}{b}\phi_{x}$ to equation $\eqref{1.1}$ we get 
\begin{equation}\label{2.9}
\left(\frac{b}{\sqrt{\xi}}u_{x}+\sqrt{\xi}\phi\right)_{x}=\frac{\r b}{\mu\sqrt{\xi}}u_{tt}+(\mu-b^{2}/\xi)\frac{\sqrt{\xi}}{\mu}\phi_{x}.
\end{equation}
Substitute $\left(\frac{b}{\sqrt{\xi}}u_{x}+\sqrt{\xi}\phi\right)_{x}$ in equation $\eqref{2.8}$ we arrive at 
\begin{equation}\label{2.10}
\begin{aligned}
&-J\int_{0}^{l}u_{ttx}\left(\frac{b}{\sqrt{\xi}}u_{x}+\sqrt{\xi}\phi\right)dx-\delta\frac{\r b}{\mu \sqrt{\xi}}\int_{0}^{l}\phi u_{ttx}dx\\
&=-\delta(\mu-b^{2}/\xi)\frac{\sqrt{\xi}}{\mu}\int_{0}^{L}|\phi_{x}|^{2}dx-\frac{\sqrt{\xi}}{2}\int_{0}^{l}\left|\frac{b}{\sqrt{\xi}}u_{x}+\sqrt{\xi}\phi\right|^{2}dx+\frac{d^{2}}{2\sqrt{\xi}}\int_{0}^{l}|w_{x}|^{2}dx.
\end{aligned}
\end{equation}
Taking into account that $u_{ttx}\phi=\frac{d}{dt}(u_{tx}\phi)-u_{tx}\phi_{t}$ and $u_{ttx}\left(\frac{b}{\sqrt{\xi}}u_{x} +\sqrt{\xi}\phi\right)=\frac{d}{dt}\left[ u_{tx}\left(\frac{b}{\sqrt{\xi}}u_{x} +\sqrt{\xi}\phi\right)\right]-u_{tx}\left(\frac{b}{\sqrt{\xi}}u_{x} +\sqrt{\xi}\phi\right)_{t}$ we obtain 
\begin{equation}\label{2.11}
\begin{aligned}
&\frac{d}{dt}\left(-J\int_{0}^{l}u_{ttx}\left(\frac{b}{\sqrt{\xi}}u_{x}+\sqrt{\xi}\phi\right)dx-\delta\frac{\r b}{\mu \sqrt{\xi}}\int_{0}^{l}\phi u_{tx}dx\right) \\
&=-\frac{Jb}{\sqrt{\xi}}\int_{0}^{l}|u_{tx}|^{2}-\underbrace{\left(J\sqrt{\xi}+\delta\frac{\r b}{\mu \sqrt{\xi}}\right)}_{C_{1}}\int_{0}^{l}\phi_{t} u_{tx}dx-\delta(\mu-b^{2}/\xi)\frac{\sqrt{\xi}}{\mu}\int_{0}^{L}|\phi_{x}|^{2}dx\\
&-\frac{\sqrt{\xi}}{2}\int_{0}^{l}\left|\frac{b}{\sqrt{\xi}}u_{x}+\sqrt{\xi}\phi\right|^{2}dx+\frac{d^{2}}{2\sqrt{\xi}}\int_{0}^{l}|w_{x}|^{2}dx.
\end{aligned}
\end{equation}
Now multiply equation $\eqref{1.3}$ by $\frac{C_{1}}{d}u_{t}$, integrate by parts over $(0,l)$ we get and using boundary conditions $\eqref{1.4}$ we have 
$$\frac{\a C_{1}}{d}\int_{0}^{l}w_{t}u_{t}dx+\frac{\kappa C_{1}}{d}\int_{0}^{l}w_{x}u_{tx}dx-C_{1}\int_{0}^{l}\phi_{t}u_{tx}dx+\frac{k C_{1}}{d}\int_{0}^{l}w u_{t}dx=0.$$
Taking into account that $w_{t}u_{t}=\frac{d}{dt}(w u_{t})-w u_{tt}$ we obtain 
\begin{equation}\label{2.12}
\frac{d}{dt}\left(\frac{\a C_{1}}{d}\int_{0}^{l}w u_{t}dx\right)=\frac{\a C_{1}}{d}\int_{0}^{l}w u_{tt}dx-\frac{\kappa C_{1}}{d}\int_{0}^{l}w_{x}u_{tx}dx+C_{1}\int_{0}^{l}\phi_{t}u_{tx}dx-\frac{k C_{1}}{d}\int_{0}^{l}w u_{t}dx.
\end{equation}
Add the two equations $\eqref{2.11}$ and $\eqref{2.12}$ and apply Young's inequality then for all $\eps_{1}$, $\eps_{2}$ and $\eps_{3} >0 $ we have 
\begin{equation}\label{2.13}
\begin{aligned}
&\frac{d}{dt}\left(-J\int_{0}^{l}u_{ttx}\left(\frac{b}{\sqrt{\xi}}u_{x}+\sqrt{\xi}\phi\right)dx-\delta\frac{\r b}{\mu \sqrt{\xi}}\int_{0}^{l}\phi u_{tx}dx+\frac{\a C_{1}}{d}\int_{0}^{l}w u_{t}dx\right)\\
&\leq-\frac{Jb}{\sqrt{\xi}}\int_{0}^{l}|u_{tx}|^{2}-\underbrace{\left(J\sqrt{\xi}+\delta\frac{\r b}{\mu \sqrt{\xi}}\right)}_{C_{1}}\int_{0}^{l}\phi_{t} u_{tx}dx-\delta(\mu-b^{2}/\xi)\frac{\sqrt{\xi}}{\mu}\int_{0}^{L}|\phi_{x}|^{2}dx-\frac{\sqrt{\xi}}{2}\int_{0}^{l}\left|\frac{b}{\sqrt{\xi}}u_{x}+\sqrt{\xi}\phi\right|^{2}dx\\
&\underbrace{\left(\frac{\a^{2}C_{1}^{2}}{2d^{2}\eps_{1}}+\frac{k^{2}C_{1}^{2}}{2d^{2}\eps_{3}}\right)}_{C_{2}}\int_{0}^{l}|w|^{2}dx+\frac{\eps_{1}}{2}\int_{0}^{l}|u_{tt}|^{2}dx+\underbrace{\left(\frac{\kappa^{2}C_{1}^{2}}{2d^{2}\eps_{2}}+\frac{d^{2}}{2\sqrt{\xi}}\right)}_{C_{3}}\int_{0}^{l}|w_{x}|^{2}dx+\frac{\eps_{2}}{2}\int_{0}^{l}|u_{tx}|^{2}dx+\frac{\eps_{3}}{2}\int_{0}^{l}|u_{t}|^{2}dx.
\end{aligned}
\end{equation}
Take $\eps_{2}=\frac{Jb}{2\sqrt{\xi}}$ we get the desired result. $\hfill{\blacksquare}$ \newline \newline

Let $$\mathcal{L}(t)=N_{1}E(t)+\mathcal{F}(t)+N_{2}\mathcal{G}(t),$$
where $N_{1}$ and $N_{2}$ are positive constants to be fixed. 
\begin{theorem}\label{thm5}
There exists positive constants $\nu_{1}$ and $\nu_{2}$ such that 
$$\nu_{1}E(t)\leq \mathcal{L}(t)\leq \nu_{2}E(t), \quad \forall t\geq 0.$$
\end{theorem}
\noindent\textbf{Proof:}
We have 
\begin{equation*}
\begin{aligned}
|\mathcal{L}(t)-N_{1}E(t)|&\leq \r\int_{0}^{l}u_{t}u dx+\frac{J\mu}{b}\int_{0}^{l}u_{xt}u_{x}dx+N_{2}J\int_{0}^{l}u_{tx}\left(\frac{b}{\sqrt{\xi}}u_{x}+\sqrt{\xi}\phi\right)dx\\
&+N_{2}\delta\frac{\r b}{\mu \sqrt{\xi}}\int_{0}^{l}u_{tx}\phi dx+N_{2}\frac{\a C_{1}}{d}\int_{0}^{l}w u_{t}dx.
\end{aligned}
\end{equation*}
Apply Young's and Poincare's inequalities we obtain 
\begin{equation*}
\begin{aligned}
|\mathcal{L}(t)-N_{1}E(t)|&\leq\left(\frac{\r}{2}+N_{2}\frac{\a C_{1}}{2d}\right)\int_{0}^{l}|u_{t}|^{2}dx+\left(\frac{\r c_{p}}{2}+\frac{J\mu}{2b}\right)\int_{0}^{l}|u_{x}|^{2}dx\\
&+\left(\frac{J\mu}{2b}+\frac{N_{2}J}{2}+N_{2}\delta \frac{\r b}{2\mu \sqrt{\xi}}\right)\int_{0}^{l}|u_{tx}|^{2}dx+\frac{N_{2}J}{2}\int_{0}^{l}\left|\frac{b}{\sqrt{\xi}}u_{x}+\sqrt{\xi}\phi\right|^{2}dx\\
&+N_{2}\delta c_{p} \frac{\r b}{2\mu \sqrt{\xi}}\int_{0}^{l}|\phi_{x}|^{2}dx+N_{2}\frac{\a C_{1}}{2d}\int_{0}^{l}|w|^{2}dx
\end{aligned}
\end{equation*}
Define 
$$N_{0} :=\max\left\{\frac{1}{\r}\left(\r+N_{3}\frac{\a C_{1}}{d}\right); \frac{1}{\mu-b^{2}/\xi}\left(\r c_{p}+\frac{J\mu}{b}\right);\right.$$
$$\left.\frac{b}{J\mu}\left(\frac{J\mu}{b}+N_{2}J+N_{2}\delta \frac{\r b}{\mu \sqrt{\xi}}\right);N_{2}J; N_{2} c_{p} \frac{\r b}{\mu \sqrt{\xi}}; N_{2}\frac{ C_{1}}{d}\right\}$$
Hence 
$$|\mathcal{L}(t)-N_{1}E(t)|\leq N_{0}E(t),$$
which implies that 
$$\nu_{1}E(t)\leq \mathcal{L}(t)\leq \nu_{2}E(t), $$ 
where  $\nu_{1}=N_{1}-N_{0}$ and $\nu_{2}=N_{1}+N_{0}$ and $N_{1}>N_{0}$. $\hfill{\blacksquare}$ \newline \newline
\noindent\textbf{Proof of Theorem \ref{Thm1}:}
It follows from Lemmas $\ref{lemma 2}$ and $\ref{lemma 3}$ that 
\begin{equation}\label{2.13}
\begin{aligned}
\frac{d}{dt}\mathcal{L}(t)&\leq -\left(\r-\frac{N_{2}\eps_{3}}{2}\right)\int_{0}^{l}|u_{t}|^{2}dx-\left(\frac{J\r}{b}-\frac{N_{2}\eps_{1}}{2}\right)\int_{0}^{l}|u_{tt}|^{2}dx-(\mu-b^{2}/\xi)\int_{0}^{l}|u_{x}|^{2}dx\\
&-\left(N_{2}\frac{Jb}{2\sqrt{\xi}}-(\frac{J\mu}{b}+2\r c_{p})\right)\int_{0}^{l}|u_{xt}|^{2}dx-\left(1+\frac{N_{2}\sqrt{\xi}}{2}\right)\int_{0}^{l}\left|\frac{b}{\sqrt{\xi}}u_{x}+\sqrt{\xi}\phi\right|^{2}dx\\
&-\left(\frac{\delta}{2}+N_{2}\delta(\mu-b^{2}/\xi)\frac{\sqrt{\xi}}{\mu}\right)\int_{0}^{l}|\phi_{x}|^{2}dx-(N_{1}k-N_{2}C_{2})\int_{0}^{l}|w|^{2}dx-(N_{1}\kappa-N_{2}C_{3})\int_{0}^{l}|w_{x}|^{2}dx
\end{aligned}
\end{equation}
Choose $\eps_{3}=\frac{\r}{N_{2}}$, $\eps_{1}=\frac{J\r}{bN_{2}}$, $N_{2}>\frac{2\sqrt{\xi}}{Jb}(\frac{J\mu}{b}+2\r c_{p})$ and $N_{1}>\max\left\{\frac{N_{2}C_{2}}{k}; \frac{N_{2}C_{3}}{\kappa}\right\}$, from where we obtain that $\zeta_{1}=\r-\frac{N_{2}\eps_{3}}{2}>0$, $\zeta_{2}=\frac{J\r}{b}-\frac{N_{2}\eps_{1}}{2}>0$, $\zeta_{3}=N_{2}\frac{Jb}{2\sqrt{\xi}}-(\frac{J\mu}{b}+2\r c_{p})>0$, $\zeta_{4}=N_{1}k-N_{2}C_{2}>0$ and $\zeta_{5}=N_{1}\kappa-N_{2}C_{3}>0$ and from where we can conclude that there exists a positive constant $\beta=2\min\{1, \zeta_{1}, \zeta_{2}, \zeta_{3}, \zeta_{4}, \zeta_{5}\}$ such that 
$$\frac{d}{dt}\mathcal{L}(t)\leq -\beta E(t),$$
by equivalence between $E(t)$ and $\mathcal{L}(t)$ according to Theorem $\ref{thm5}$ we get:
$$ \frac{d}{dt}\mathcal{L}(t)\leq -\omega \mathcal{L}(t),$$
where $\omega=\frac{\beta}{\nu_{2}}$.
Now integrate the above inequality over $(0,t)$ we obtain
$$ \mathcal{L}(t)\leq \mathcal{L}(0)e^{-\omega t},$$
again by equivalence between $E(t)$ and $\mathcal{L}(t)$ according to Theorem $\ref{thm5}$ we arrive at 
$$E(t)\leq ME(0)e^{-\omega t},$$
where $M=\frac{\nu_{2}}{\nu_{1}}$. $\hfill{\blacksquare}$ 

\section{Numerical approximation}
\SE
First we denote by $\hat{u}=u_{t},\hat{\phi}=\phi_{t}$, $\hat{w}=w_{t}$  and we introduce the following weak form after multiplying the equations $\eqref{1.1}$ and $\eqref{1.2}$ by $\bar{u}, \bar{\phi}, \bar{w} \in H_{0}^{1}(0,l)$ 
\begin{equation}\label{WP}
(WP)\left\{
\begin{aligned}
&\r(\hat{u}_{t},\bar{u})+\mu(u_{x},\bar{u}_{x})+b(\phi,\bar{u}_{x})=0\\
&J(\hat{u}_{t},\bar{\phi}_{x})+\delta(\phi_{x},\bar{\phi}_{x})+b(u_{x},\bar{\phi})+\xi(\phi, \bar{\phi})+d(w_{x},\bar{\phi}) =0\\
&\a(\hat{w}, \bar{w})+\kappa(w_{x}, \bar{w}_{x})+d(\hat{\phi}_{x}, \bar{w})+k(w,\bar{w})=0
\end{aligned}
\right.
\end{equation}
Let us partition the interval $(0;l)$ into subintervals $I_{j} = (x_{j-1}; x_{j})$ of length $h = \frac{1}
{s}$ with $0 = x_{0} < x_{1} < ... < x_{s} = l$ and define the associated finite element spaces by 
$$S_{h}^{0}=\{u\in H_{0}^{1}(0,l); u\in C([0,l]), u|_{I_{j}}\in P_{1}(K)\}.$$  
For a given final time $T$ and a positive integer $N$, define the time step $\Delta t=\frac{T}{N}$ and the nodes $t_{n}=n\Delta t, n=0,...,N$. By using the Implicit Euler scheme in time and the finite element variational approximation in space, we introduce the following scheme. For $\bar{u}_{h}, \bar{\phi}_{h}, \bar{w}_{h}\in S_{h}^{0}$  ,  find $u_{h}^{n}, \phi_{h}^{n}, w_{h}^{n} \in S_{h}^{0}$  such that 
\begin{equation}\label{NP}
(NP)\left\{
\begin{aligned}
&\frac{\r}{\Delta t}(\hat{u}_{h}^{n}-\hat{u}_{h}^{n-1},\bar{u}_{h})+\mu(u_{hx}^{n},\bar{u}_{hx})+b(\phi_{h}^{n},\bar{u}_{hx})=0\\
&\frac{J}{\Delta t}(\hat{u}_{h}^{n}-\hat{u}_{h}^{n-1},\bar{\phi}_{hx})+\delta(\phi_{hx}^{n},\bar{\phi}_{hx})+b(u_{hx}^{n},\bar{\phi}_{h})+\xi(\phi_{h}^{n},\bar{\phi}_{h})+d(w_{hx}^{n},\bar{\phi}_{h})=0\\
&\frac{\a}{\Delta t}(w_{h}^{n}-w_{h}^{n-1}, \bar{w}_{h})+\kappa(w_{hx}^{n}, \bar{w}_{hx})+d(\hat{\phi}_{hx}^{n}, \bar{w}_{h})+k(w_{h}^{n},\bar{w}_{h})=0
\end{aligned}
\right.
\end{equation}
Where $u_{h}^{n}=u_{h}^{n-1}+\Delta t \hat{u}_{h}^{n}$;  $\phi_{h}^{n}=\phi_{h}^{n-1}+\Delta t\hat{\phi}_{h}^{n}$ and $w_{h}^{n}=w_{h}^{n-1}+\Delta t \hat{w}_{h}^{n}$.\newline For a continuous function $f(t)$, let $f^{n}=f(t_{n})$ and for a sequence $\{f_{n}\}_{n=1}^{N}$ let $\sigma f^{n}=(f^{n}-f^{n-1})/
\Delta t$ \newline 
The discrete energy at certain time $t_{n}$ is defined by 
$$E_{n}=\frac{1}{2}\left(\r||\hat{u}_{h}^{n}||^{2}+\frac{J\r}{b}||\sigma\hat{u}_{h}^{n}||^{2}+(\mu-b^{2}/\xi)||u_{hx}^{n}||^{2}+\frac{J\mu}{b}||\hat{u}_{hx}^{n}||^{2}+\delta||\phi_{hx}^{n}||^{2}+||\frac{b}{\sqrt{\xi}}u_{hx}^{n}+\sqrt{\xi}\phi_{h}^{n}||^{2}+\a||w_{h}^{n}||^{2}\right)$$
Where $||.||$ denotes the $L^{2}$--norm. The decay of energy is presented in the following proposition.
\begin{proposition}
For all $n=0,...,N$, we have 
$$\frac{E_{n}-E_{n-1}}{\Delta t}<0$$
\end{proposition}

\noindent{\textbf{Proof}} Taking $\bar{u}_{h}=\hat{u}_{h}^{n}$,  $\bar{\phi}_{h}=\hat{\phi}_{h}^{n}$ and $\bar{w}_{h}=w_{h}^{n}$ in $(NP)$ we get 
\begin{equation}\label{n2.3}
\left\{
\begin{aligned}
&\frac{\r}{2\Delta t}\left(||\hat{u}_{h}^{n}-\hat{u}_{h}^{n-1}||^{2}+||\hat{u}_{h}^{n}||^{2}-||\hat{u}_{h}^{n-1}||^{2}\right)+\frac{\mu}{2\Delta t}\left(||u_{hx}^{n}-u_{hx}^{n-1}||^{2}+||u_{hx}^{n}||^{2}-||u_{hx}^{n-1}||^{2}\right)\\
&+b(\phi_{h}^{n},\hat{u}_{hx}^{n})=0\\
&\frac{J}{\Delta t}\left(\hat{u}_{h}^{n}-\hat{u}_{h}^{n-1},\hat{\phi}_{hx}^{n}\right)+\frac{\delta}{2\Delta t}\left(||\phi_{hx}^{n}-\phi_{hx}^{n-1}||^{2}+||\phi_{hx}^{n}||^{2}-||\phi_{hx}^{n-1}||^{2}\right)\\
&+b(u_{hx}^{n},\hat{\phi}_{h}^{n})+\frac{\xi}{2\Delta t}\left(||\phi_{h}^{n}-\phi_{h}^{n-1}||^{2}+||\phi_{h}^{n}||^{2}-||\phi_{h}^{n-1}||^{2}\right)+d(w_{hx}^{n},\hat{\phi}_{h}^{n})=0\\
&\frac{\a}{2\Delta t}\left(||w_{h}^{n}-w_{h}^{n-1}||^{2}+||w_{h}^{n}||^{2}-||w_{h}^{n-1}||^{2}\right)+\kappa||w_{hx}^{n}||^{2}+d(\hat{\phi}_{hx}^{n}, w_{h}^{n})+k||w_{h}^{n}||^{2}=0
\end{aligned}
\right.
\end{equation}
Sum the above two equations and note that
\begin{equation*}
\begin{aligned} 
&b(\phi_{h}^{n},\hat{u}_{hx}^{n})+b(u_{hx}^{n},\hat{\phi}_{h}^{n})\\
&=\frac{b}{\Delta t}\left( (\phi_{h}^{n}, u_{hx}^{n}-u_{hx}^{n-1})+ (u_{hx}^{n},\phi_{h}^{n}-\phi_{h}^{n-1})\right)\\
&=\frac{b}{\Delta t}\left( (\phi_{h}^{n}, u_{hx}^{n})-(\phi_{h}^{n-1}, u_{hx}^{n-1})+(\phi_{h}^{n}-\phi_{h}^{n-1}, u_{hx}^{n}-u_{hx}^{n-1})\right)
\end{aligned}
\end{equation*}
Now by the inequality $\mu\xi-b^{2}>0$ we get 
\begin{equation*}
\begin{aligned}
&\frac{\mu}{2\Delta t}||u_{hx}^{n}-u_{hx}^{n-1}||^{2}+\frac{\xi}{2\Delta t}||\phi_{h}^{n}-\phi_{h}^{n-1}||^{2}+\frac{b}{\Delta t}(\phi_{h}^{n}-\phi_{h}^{n-1}, u_{hx}^{n}-u_{hx}^{n-1})\\
>&\frac{1}{2 \Delta t}\left( \frac{b^{2}}{\xi}||u_{hx}^{n}-u_{hx}^{n-1}||^{2}+\xi||\phi_{h}^{n}-\phi_{h}^{n-1}||^{2}+2b(\phi_{h}^{n}-\phi_{h}^{n-1}, u_{hx}^{n}-u_{hx}^{n-1})\right)\\
=&\frac{1}{2\Delta t}||\frac{b}{\sqrt{\xi}}(u_{hx}^{n}-u_{hx}^{n-1})+\sqrt{\xi}(\phi_{h}^{n}-\phi_{h}^{n-1})||^{2}>0\\
\end{aligned}
\end{equation*}
And by the positivity of the terms $||\hat{u}_{h}^{n}-\hat{u}_{h}^{n-1}||^{2}$, $||\phi_{hx}^{n}-\phi_{hx}^{n-1}||^{2}$, $||w_{h}^{n}-w_{h}^{n-1}||^{2}$, $||w_{hx}^{n}||^{2}$ and $||w_{h}^{n}||^{2}$   we get the following inequality
\begin{equation}\label{n2.4}
\begin{aligned} 
&\frac{\r}{2\Delta t}\left(||\hat{u}_{h}^{n}||^{2}-||\hat{u}_{h}^{n-1}||^{2}\right)+ \frac{\mu}{2\Delta t}\left(||u_{hx}^{n}||^{2}-||u_{hx}^{n-1}||^{2}\right)+\frac{J}{\Delta t}\left(\hat{u}_{h}^{n}-\hat{u}_{h}^{n-1},\hat{\phi}_{hx}^{n}\right)+\frac{\delta}{2\Delta t}\left(||\phi_{hx}^{n}||^{2}-||\phi_{hx}^{n-1}||^{2}\right)\\
&+\frac{\xi}{2\Delta t}\left(||\phi_{h}^{n}||^{2}-||\phi_{h}^{n-1}||^{2}\right)+\frac{\a}{2\Delta t}\left(||w_{h}^{n}||^{2}-||w_{h}^{n-1}||^{2}\right)+\frac{b}{\Delta t}\left( (\phi_{h}^{n}, u_{hx}^{n})-(\phi_{h}^{n-1}, u_{hx}^{n-1})\right)<0.
\end{aligned}
\end{equation}   
Add then subtract $\frac{b^{2}}{2\xi\Delta t}(||u_{hx}^{n}||^{2}-||u_{hx}^{n-1}||^{2})$ to the above inequality we get: 
\begin{equation}\label{n2.5}
\begin{aligned} 
&\frac{\r}{2\Delta t}\left(||\hat{u}_{h}^{n}||^{2}-||\hat{u}_{h}^{n-1}||^{2}\right)+ \frac{(\mu-b^{2}/\xi)}{2\Delta t}\left(||u_{hx}^{n}||^{2}-||u_{hx}^{n-1}||^{2}\right)+\frac{J}{\Delta t}\left(\hat{u}_{h}^{n}-\hat{u}_{h}^{n-1},\hat{\phi}_{hx}^{n}\right)+\frac{\delta}{2\Delta t}\left(||\phi_{hx}^{n}||^{2}-||\phi_{hx}^{n-1}||^{2}\right)\\
&+\frac{1}{2\Delta t}\left(||\frac{b}{\sqrt{\xi}}u_{hx}^{n}+\sqrt{\xi}\phi_{h}^{n}||^{2}-||\frac{b}{\sqrt{\xi}}u_{hx}^{n-1}+\sqrt{\xi}\phi_{h}^{n-1}||^{2}\right)+\frac{\a}{2\Delta t}\left(||w_{h}^{n}||^{2}-||w_{h}^{n-1}||^{2}\right)<0.
\end{aligned}
\end{equation}   
Moreover by equation $\eqref{1.1}$ we get 
$$\phi_{xt}=\frac{\r}{b}u_{ttt}-\frac{\mu}{b}u_{xxt},$$
which implies that 
\begin{equation}\label{n2.6}
\begin{aligned}
&\frac{J}{\Delta t}\left(\hat{u}_{h}^{n}-\hat{u}_{h}^{n-1},\hat{\phi}_{hx}^{n}\right)\\
=&\frac{J}{\Delta t}\left(\hat{u}_{h}^{n}-\hat{u}_{h}^{n-1},\frac{\r}{b}\sigma^{2}\hat{u}_{h}^{n}-\frac{\mu}{b}\hat{u}_{hxx}^{n}\right)\\
=&\frac{J\r}{b\Delta t^{2}}\left(\hat{u}_{h}^{n}-\hat{u}_{h}^{n-1}, \sigma\hat{u}_{h}^{n}-\sigma\hat{u}_{h}^{n-1}\right)
+\frac{J\mu}{b\Delta t}\left(\hat{u}_{hx}^{n}-\hat{u}_{hx}^{n-1},\hat{u}_{hx}^{n}\right)\\
=& \frac{J\r}{b\Delta t^{3}}\left(\hat{u}_{h}^{n}-\hat{u}_{h}^{n-1}, \hat{u}_{h}^{n}-\hat{u}_{h}^{n-1}-(\hat{u}_{h}^{n-1}-\hat{u}_{h}^{n-2})\right)
+\frac{J\mu}{2b\Delta t}\left( ||\hat{u}_{hx}^{n}-\hat{u}_{hx}^{n-1}||^{2}+||\hat{u}_{hx}^{n}||^{2}-||\hat{u}_{hx}^{n-1}||^{2}\right)\\
= &\frac{J\r}{2b\Delta t^{3}}\left( ||\hat{u}_{h}^{n}-\hat{u}_{h}^{n-1}-(\hat{u}_{h}^{n-1}-\hat{u}_{h}^{n-2})||^{2}+||\hat{u}_{h}^{n}-\hat{u}_{h}^{n-1}||^{2}-||\hat{u}_{h}^{n-1}-\hat{u}_{h}^{n-2}||^{2}\right)\\
&+\frac{J\mu}{2b\Delta t}\left( ||\hat{u}_{hx}^{n}-\hat{u}_{hx}^{n-1}||^{2}+||\hat{u}_{hx}^{n}||^{2}-||\hat{u}_{hx}^{n-1}||^{2}\right)\\
\geq &\frac{J\r}{2b\Delta t^{3}}\left(||\hat{u}_{h}^{n}-\hat{u}_{h}^{n-1}||^{2}-||\hat{u}_{h}^{n-1}-\hat{u}_{h}^{n-2}||^{2}\right)
+\frac{J\mu}{2b\Delta t}\left(||\hat{u}_{hx}^{n}||^{2}-||\hat{u}_{hx}^{n-1}||^{2}\right)
\end{aligned}
\end{equation}
Now substitute inequality $\eqref{n2.6}$ in equation $\eqref{n2.5}$ we get 
\begin{equation}\label{n2.7}
\begin{aligned} 
&\frac{\r}{2\Delta t}\left(||\hat{u}_{h}^{n}||^{2}-||\hat{u}_{h}^{n-1}||^{2}\right)+ \frac{(\mu-b^{2}/\xi)}{2\Delta t}\left(||u_{hx}^{n}||^{2}-||u_{hx}^{n-1}||^{2}\right)+\frac{J\r}{2b\Delta t^{3}}\left(||\hat{u}_{h}^{n}-\hat{u}_{h}^{n-1}||^{2}-||\hat{u}_{h}^{n-1}-\hat{u}_{h}^{n-2}||^{2}\right)\\
&+\frac{J\mu}{2b\Delta t}\left(||\hat{u}_{hx}^{n}||^{2}-||\hat{u}_{hx}^{n-1}||^{2}\right)+\frac{\delta}{2\Delta t}\left(||\phi_{hx}^{n}||^{2}-||\phi_{hx}^{n-1}||^{2}\right)+\frac{1}{2\Delta t}\left(||\frac{b}{\sqrt{\xi}}u_{hx}^{n}+\sqrt{\xi}\phi_{h}^{n}||^{2}-||\frac{b}{\sqrt{\xi}}u_{hx}^{n-1}+\sqrt{\xi}\phi_{h}^{n-1}||^{2}\right)\\
&+\frac{\a}{2\Delta t}\left(||w_{h}^{n}||^{2}-||w_{h}^{n-1}||^{2}\right)<0.
\end{aligned}
\end{equation}   
which is the desired result. $\hfill{\blacksquare}$

\section{A Priori Error estimate}
\SE
\begin{theorem}
Suppose that the solution $(u,\phi)$ of equations $\eqref{1.1}$--$\eqref{1.4}$ belongs to the space 
$$(H^{4}(0,T;H^{2}(0,l)))^{2}$$
then for $\Delta t$ sufficiently small, the following a priori error estimate holds 
$$||\hat{u}_{h}^{n}-u_{t}(t_{n})||^{2}+||\hat{\phi}_{h}^{n}-\phi_{t}(t_{n})||^{2}+||u_{hx}^{n}-u_{x}(t_{n})||^{2}+||\phi_{hx}^{n}-\phi_{x}(t_{n})||^{2}+||\phi_{h}^{n}-\phi(t_{n})||^{2}+||w_{h}^{n}-w(t_{n})||^{2}\leq c(h^{2}+(\Delta t)^{2})$$
\end{theorem}

\noindent{\textbf{Proof:}} We start by introducing the following terms 
$$e^{n}=u_{h}^{n}-P_{h}^{0}u(t_{n}),$$ 
$$\hat{e}^{n}=\hat{u}_{h}^{n}-P_{h}^{0}u_{t}(t_{n}),$$
$$q^{n}=\phi_{h}^{n}-P_{h}^{0}\phi(t_{n}),$$ 
$$ \hat{q}^{n}=\hat{\phi}_{h}^{n}-P_{h}^{0}\phi_{t}(t_{n}),$$
$$R^{n}=w_{h}^{n}-P_{h}^{0}w(t_{n}).$$
Where $P_{h}^{0}$ is projection operator such that 
$$P_{h}^{0}:H_{0}^{1}(0,l)\longrightarrow S_{h}^{0},$$
defined by 
$$\left( (P_{h}^{0}u)_{x}, \chi_{x}\right)=\left(u_{x}, \chi_{x}\right),\quad \mbox{forall} \quad \chi \in S_{h}^{0}, u \in H_{0}^{1}(0,l) $$
satisfying 
$$|u- P_{h}^{0}u|\leq Ch|u_{x}|, \quad \mbox{and} \quad (P_{h}^{0}u)(x_{i})=u(x_{i}).$$

\textbf{Step 1.} Replace $u_{h}^{n}$ by $e^{n}+P_{h}^{0}u(t_{n})$, $\hat{u}_{h}^{n}$ by $\hat{e}^{n}+P_{h}^{0}u_{t}(t_{n})$, $\phi_{h}^{n}$ by $q^{n}+P_{h}^{0}\phi(t_{n})$ and take $\bar{u}_{h}=\hat{e}^{n}$ in the first equation of $(NP)$ we get 
\begin{equation}\label{3.1}
\begin{aligned}
&\frac{\r}{\Delta t}(\hat{e}^{n}-\hat{e}^{n-1},\hat{e}^{n})+\frac{\r}{\Delta t}(P_{h}^{0}u_{t}(t_{n})-P_{h}^{0}u_{t}(t_{n-1}),\hat{e}^{n})+\mu(e_{x}^{n},\hat{e}_{x}^{n})+\mu((P_{h}^{0}u(t_{n}))_{x},\hat{e}_{x}^{n})+b(q^{n},\hat{e}_{x}^{n})+b(P_{h}^{0}\phi(t_{n}),\hat{e}_{x}^{n})=0,
\end{aligned}
\end{equation}
which implies
\begin{equation}\label{3.2}
\begin{aligned}
&\frac{\r}{2\Delta t}\left(||\hat{e}^{n}-\hat{e}^{n-1}||^{2}+||\hat{e}^{n}||^{2}-||\hat{e}^{n-1}||^{2}\right)+\frac{\r}{\Delta t}\left(P_{h}^{0}u_{t}(t_{n})-P_{h}^{0}u_{t}(t_{n-1}),\hat{e}^{n}\right)+\mu\left(e_{x}^{n},\hat{e}_{x}^{n}\right)+\mu\left((P_{h}^{0}u(t_{n}))_{x},\hat{e}_{x}^{n}\right)\\
&+b\left(q^{n},\hat{e}_{x}^{n}\right)+b\left(P_{h}^{0}\phi(t_{n}),\hat{e}_{x}^{n}\right)=0.
\end{aligned}
\end{equation}
 Now replace $\bar{u}$ by $\hat{e}^{n}$ in the first equation of $(WP)$ at time $t_{n}$ and then combine it with equation $\eqref{3.2}$ we obtain 
\begin{equation}\label{3.3}
\begin{aligned}
&\frac{\r}{2\Delta t}\left(||\hat{e}^{n}-\hat{e}^{n-1}||^{2}+||\hat{e}^{n}||^{2}-||\hat{e}^{n-1}||^{2}\right)+\mu\left(e_{x}^{n},\hat{e}_{x}^{n}\right)+b\left(q^{n},\hat{e}_{x}^{n}\right)\\
&=\r\left(u_{tt}(t_{n})-\frac{P_{h}^{0}u_{t}(t_{n})-P_{h}^{0}u_{t}(t_{n-1})}{\Delta t},\hat{e}^{n}\right)+\mu\left(u_{x}(t_{n})-(P_{h}^{0}u(t_{n}))_{x},\hat{e}_{x}^{n}\right)+b\left(\phi(t_{n})-P_{h}^{0}\phi(t_{n}),\hat{e}_{x}^{n}\right)
\end{aligned}
\end{equation}
\textbf{Step 2.}  Apply the same procedure to the second equation of $(NP)$ that is replace $\phi_{h}^{n}$ by $q^{n}+P_{h}^{0}\phi(t_{n})$, $\hat{\phi}_{h}^{n}$ by $\hat{q}^{n}+P_{h}^{0}\phi_{t}(t_{n})$, $\hat{u}_{h}^{n}$ by $\hat{e}^{n}+P_{h}^{0}u_{t}(t_{n})$ and $w_{h}^{n}$ by $R^{n}+P_{h}^{0}w(t_{n})$ and for $\bar{\phi}_{h}=\hat{q}^{n}$ we have 
\begin{equation}\label{3.4}
\begin{aligned}
&\frac{J}{\Delta t}\left(\hat{e}^{n}-\hat{e}^{n-1},\hat{q}_{x}^{n}\right)+\frac{J}{\Delta t}(P_{h}^{0}u_{t}(t_{n})-P_{h}^{0}u_{t}(t_{n-1}),\hat{q}_{x}^{n})+\delta(q_{x}^{n},\hat{q}_{x}^{n})+\delta\left((P_{h}^{0}\phi(t_{n}))_{x},\hat{q}_{x}^{n}\right)+b(e_{x}^{n},\hat{q}^{n})+b\left((P_{h}^{0}u(t_{n}))_{x},\hat{q}^{n}\right)\\
&+\xi(q^{n},\hat{q}^{n})+\xi(P_{h}^{0}\phi(t_{n}),\hat{q}^{n})+ d(R_{x}^{n},\hat{q}^{n})+ d\left((P_{h}^{0}w(t_{n}))_{x},\hat{q}^{n}\right)=0
\end{aligned}
\end{equation}
Taking into consideration that 
$$\hat{q}_{x}^{n}=\frac{\r}{b}\sigma^{2}\hat{e}^{n}-\frac{\mu}{b}\hat{e}_{xx}^{n},$$
then 
\begin{equation*}
\begin{aligned}
&\frac{J}{\Delta t}\left(\hat{e}^{n}-\hat{e}^{n-1},\hat{q}_{x}^{n}\right)\\
=& \frac{J}{\Delta t}\left(\hat{e}^{n}-\hat{e}^{n-1},\frac{\r}{b}\sigma^{2}\hat{e}^{n}-\frac{\mu}{b}\hat{e}_{xx}^{n}\right)\\
=&\frac{J\r}{b\Delta t^{2}}\left(\hat{e}^{n}-\hat{e}^{n-1},\sigma\hat{e}^{n}-\sigma\hat{e}^{n-1}\right)+\frac{J\mu}{b\Delta t}\left(\hat{e}_{x}^{n}-\hat{e}_{x}^{n-1},\hat{e}_{x}^{n}\right)\\
=&\frac{J\r}{b\Delta t^{3}}\left(\hat{e}^{n}-\hat{e}^{n-1},\hat{e}^{n}-\hat{e}^{n-1}-(\hat{e}^{n-1}-\hat{e}^{n-2})\right)+\frac{J\mu}{b\Delta t}\left(\hat{e}_{x}^{n}-\hat{e}_{x}^{n-1},\hat{e}_{x}^{n}\right)\\
=&\frac{J\r}{2b\Delta t^{3}}\left(||\hat{e}^{n}-\hat{e}^{n-1}||^{2}-||\hat{e}^{n-1}-\hat{e}^{n-2}||^{2}+||\hat{e}^{n}-\hat{e}^{n-1}-(\hat{e}^{n-1}-\hat{e}^{n-2})||^{2}\right)+\frac{J\mu}{2b\Delta t}\left(||\hat{e}_{x}^{n}||^{2}-||\hat{e}_{x}^{n-1}||^{2}+||\hat{e}_{x}^{n}-\hat{e}_{x}^{n-1}||^{2}\right)\\
=&\frac{J\r}{2b\Delta t}\left(||\sigma\hat{e}^{n}||^{2}-||\sigma\hat{e}^{n-1}||^{2}+||\sigma(\hat{e}^{n}-\hat{e}^{n-1})||^{2}\right)+\frac{J\mu}{2b\Delta t}\left(||\hat{e}_{x}^{n}||^{2}-||\hat{e}_{x}^{n-1}||^{2}+||\hat{e}_{x}^{n}-\hat{e}_{x}^{n-1}||^{2}\right)
\end{aligned}
\end{equation*}
and 
\begin{equation*}
\begin{aligned}
&\frac{J}{\Delta t}(P_{h}^{0}u_{t}(t_{n})-P_{h}^{0}u_{t}(t_{n-1}),\hat{q}_{x}^{n})\\
=&\frac{J}{\Delta t}(P_{h}^{0}u_{t}(t_{n})-P_{h}^{0}u_{t}(t_{n-1}),\frac{\r}{b}\sigma^{2}\hat{e}^{n}-\frac{\mu}{b}\hat{e}_{xx}^{n})\\
=&\frac{J\r}{b\Delta t}(P_{h}^{0}u_{t}(t_{n})-P_{h}^{0}u_{t}(t_{n-1}),\sigma^{2}\hat{e}^{n})+\frac{J\mu}{b\Delta t}\left((P_{h}^{0}u_{t}(t_{n}))_{x}-(P_{h}^{0}u_{t}(t_{n-1}))_{x},\hat{e}_{x}^{n}\right)
\end{aligned}
\end{equation*}
Then $\eqref{3.4}$ becomes 
\begin{equation}\label{3.4*}
\begin{aligned}
&\frac{J\r}{2b\Delta t}\left(||\sigma\hat{e}^{n}||^{2}-||\sigma\hat{e}^{n-1}||^{2}+||\sigma(\hat{e}^{n}-\hat{e}^{n-1})||^{2}\right)+\frac{J\mu}{2b\Delta t}\left(||\hat{e}_{x}^{n}||^{2}-||\hat{e}_{x}^{n-1}||^{2}+||\hat{e}_{x}^{n}-\hat{e}_{x}^{n-1}||^{2}\right)\\
&+\frac{J\r}{b\Delta t}(P_{h}^{0}u_{t}(t_{n})-P_{h}^{0}u_{t}(t_{n-1}),\sigma^{2}\hat{e}^{n})+\frac{J\mu}{b\Delta t}\left((P_{h}^{0}u_{t}(t_{n}))_{x}-(P_{h}^{0}u_{t}(t_{n-1}))_{x},\hat{e}_{x}^{n}\right)+\delta(q_{x}^{n},\hat{q}_{x}^{n})+\delta\left((P_{h}^{0}\phi(t_{n}))_{x},\hat{q}_{x}^{n}\right)\\
&+b(e_{x}^{n},\hat{q}^{n})+b\left((P_{h}^{0}u(t_{n}))_{x},\hat{q}^{n}\right)+\xi(q^{n},\hat{q}^{n})+\xi(P_{h}^{0}\phi(t_{n}),\hat{q}^{n})+ d(R_{x}^{n},\hat{q}^{n})+ d\left((P_{h}^{0}w(t_{n}))_{x},\hat{q}^{n}\right)=0
\end{aligned}
\end{equation}

 Now replace $\bar{\phi}$ by $\hat{q}^{n}$ in the second equation of $(WP)$ at time $t_{n}$ and using $\hat{q}_{x}^{n}=\frac{\r}{b}\sigma^{2}\hat{e}^{n}-\frac{\mu}{b}\hat{e}_{xx}^{n}$ we get 
 \begin{equation}\label{3.4**}
 \frac{J\r}{b}\left(\hat{u}_{t},\sigma^{2}\hat{e}^{n}\right)+\frac{J\mu}{b}\left(\hat{u}_{tx},\hat{e}_{x}^{n}\right)+\delta(\phi_{x},\hat{q}_{x}^{n})+b(u_{x},\hat{q}^{n})+\xi(\phi, \hat{q}^{n})+d(w_{x},\hat{q}^{n}) =0.
\end{equation}
 
 Combine $\eqref{3.4**}$ with equation $\eqref{3.4*}$ we obtain 
 \begin{equation}\label{3.4***}
\begin{aligned}
&\frac{J\r}{2b\Delta t}\left(||\sigma\hat{e}^{n}||^{2}-||\sigma\hat{e}^{n-1}||^{2}+||\sigma(\hat{e}^{n}-\hat{e}^{n-1})||^{2}\right)+\frac{J\mu}{2b\Delta t}\left(||\hat{e}_{x}^{n}||^{2}-||\hat{e}_{x}^{n-1}||^{2}+||\hat{e}_{x}^{n}-\hat{e}_{x}^{n-1}||^{2}\right)\\
&+\delta(q_{x}^{n},\hat{q}_{x}^{n})+b(e_{x}^{n},\hat{q}^{n})+\xi(q^{n},\hat{q}^{n})+ d(R_{x}^{n},\hat{q}^{n})\\
&=\frac{J\r}{b}\left(u_{tt}(t_{n})-\frac{P_{h}^{0}u_{t}(t_{n})-P_{h}^{0}u_{t}(t_{n-1})}{\Delta t},\sigma^{2}\hat{e}^{n}\right)+\frac{J\mu}{b}\left(u_{ttx}(t_{n})-\frac{(P_{h}^{0}u_{t}(t_{n}))_{x}-(P_{h}^{0}u_{t}(t_{n-1}))_{x}}{\Delta t},\hat{e}_{x}^{n}\right)\\
&+\delta\left(\phi_{x}(t_{n})-(P_{h}^{0}\phi(t_{n}))_{x},\hat{q}_{x}^{n}\right)+b\left(u_{x}(t_{n})-(P_{h}^{0}u(t_{n}))_{x},\hat{q}^{n}\right)+\xi(\phi(t_{n})-P_{h}^{0}\phi(t_{n}),\hat{q}^{n})+d\left(w_{x}(t_{n})-(P_{h}^{0}w(t_{n}))_{x},\hat{q}^{n}\right)
\end{aligned}
\end{equation}
Now in the third equation of $(NP)$  replace $\phi_{h}^{n}$ by $q^{n}+P_{h}^{0}\phi(t_{n})$, $\hat{\phi}_{h}^{n}$ by $\hat{q}^{n}+P_{h}^{0}\phi_{t}(t_{n})$ and $w_{h}^{n}$ by $R^{n}+P_{h}^{0}w(t_{n})$ and for $\bar{w}_{h}=R^{n}$ we have 
\begin{equation}\label{3.44}
\begin{aligned}
& \frac{\a}{\Delta t}(R^{n}-R^{n-1}, R^{n})+\frac{\a}{\Delta t}(P_{h}^{0}w(t_{n})-P_{h}^{0}w(t_{n-1}), R^{n})+\kappa(R_{x}^{n}, R_{x}^{n})+\kappa\left((P_{h}^{0}w(t_{n}))_{x}, R_{x}^{n}\right)\\
&+d(\hat{q}_{x}^{n}, R^{n})+d\left((P_{h}^{0}\phi_{t}(t_{n}))_{x}, R^{n}\right)+k(R^{n},R^{n})+k(P_{h}^{0}w(t_{n}),R^{n})=0
\end{aligned}
 \end{equation}
 Now replace $\bar{w}$ by $R^{n}$ in the third equation of $(WP)$ at time $t_{n}$ and then combine it with equation $\eqref{3.44}$ we obtain
 \begin{equation}\label{3.45}
 \begin{aligned}
 & \frac{\a}{2\Delta t}\left(||R^{n}-R^{n-1}||^{2}+|| R^{n}||^{2}-|| R^{n-1}||^{2}\right)+\kappa||R_{x}^{n}||^{2}+d(\hat{q}_{x}^{n}, R^{n})+k||R^{n}||^{2}\\
&=\a\left(w_{t}(t_{n})-\frac{P_{h}^{0}w(t_{n})-P_{h}^{0}w(t_{n-1})}{\Delta t}, R^{n}\right)+\kappa\left(w_{x}(t_{n})-(P_{h}^{0}w(t_{n}))_{x}, R_{x}^{n}\right)+d\left(\phi_{tx}(t_{n})-(P_{h}^{0}\phi_{t}(t_{n}))_{x}, R^{n}\right)\\
&+k(w-P_{h}^{0}w(t_{n}),R^{n})
\end{aligned}
\end{equation}

\textbf{Step 3.}
Summing equations $\eqref{3.3}$ and $\eqref{3.4***}$ we obtain 
\begin{equation}\label{3.5}
\begin{aligned}
&\frac{\r}{2\Delta t}\left(||\hat{e}^{n}-\hat{e}^{n-1}||^{2}+||\hat{e}^{n}||^{2}-||\hat{e}^{n-1}||^{2}\right)+\frac{J\r}{2b\Delta t}\left(||\sigma\hat{e}^{n}||^{2}-||\sigma\hat{e}^{n-1}||^{2}+||\sigma(\hat{e}^{n}-\hat{e}^{n-1})||^{2}\right)\\
&+\frac{J\mu}{2b\Delta t}\left(||\hat{e}_{x}^{n}||^{2}-||\hat{e}_{x}^{n-1}||^{2}+||\hat{e}_{x}^{n}-\hat{e}_{x}^{n-1}||^{2}\right)+\frac{\a}{2\Delta t}\left(||R^{n}-R^{n-1}||^{2}+|| R^{n}||^{2}-|| R^{n-1}||^{2}\right)\\
&+\mu\left(e_{x}^{n},\hat{e}_{x}^{n}\right)+\delta(q_{x}^{n},\hat{q}_{x}^{n})+b\left(q^{n},\hat{e}_{x}^{n}\right)+b(e_{x}^{n},\hat{q}^{n})+\xi(q^{n},\hat{q}^{n})+\kappa||R_{x}^{n}||^{2}+k||R^{n}||^{2}\\
&=\r\left(u_{tt}(t_{n})-\frac{P_{h}^{0}u_{t}(t_{n})-P_{h}^{0}u_{t}(t_{n-1})}{\Delta t},\hat{e}^{n}\right)+\frac{J\r}{b}\left(u_{tt}(t_{n})-\frac{P_{h}^{0}u_{t}(t_{n})-P_{h}^{0}u_{t}(t_{n-1})}{\Delta t},\sigma^{2}\hat{e}^{n}\right)\\
&+\frac{J\mu}{b}\left(u_{ttx}(t_{n})-\frac{(P_{h}^{0}u_{t}(t_{n}))_{x}-(P_{h}^{0}u_{t}(t_{n-1}))_{x}}{\Delta t},\hat{e}_{x}^{n}\right)+\a\left(w_{t}(t_{n})-\frac{P_{h}^{0}w(t_{n})-P_{h}^{0}w(t_{n-1})}{\Delta t}, R^{n}\right)\\
&+\mu\left(u_{x}(t_{n})-(P_{h}^{0}u(t_{n}))_{x},\hat{e}_{x}^{n}\right)+\delta\left(\phi_{x}(t_{n})-(P_{h}^{0}\phi(t_{n}))_{x},\hat{q}_{x}^{n}\right)+b\left(\phi(t_{n})-P_{h}^{0}\phi(t_{n}),\hat{e}_{x}^{n}\right)+b\left(u_{x}(t_{n})-(P_{h}^{0}u(t_{n}))_{x},\hat{q}^{n}\right)\\
&+\xi(\phi(t_{n})-P_{h}^{0}\phi(t_{n}),\hat{q}^{n})+d\left(w_{x}(t_{n})-(P_{h}^{0}w(t_{n}))_{x},\hat{q}^{n}\right)+\kappa\left(w_{x}(t_{n})-(P_{h}^{0}w(t_{n}))_{x}, R_{x}^{n}\right)\\
&+d\left(\phi_{tx}(t_{n})-(P_{h}^{0}\phi_{t}(t_{n}))_{x}, R^{n}\right)+k(w-P_{h}^{0}w(t_{n}),R^{n}).
\end{aligned}
\end{equation}
\textbf{Step 4.}
Note that 
\begin{equation}\label{3.6}
\begin{aligned}
&b\left(q^{n},\hat{e}_{x}^{n}\right)+b\left(e_{x}^{n},\hat{q}^{n}\right)\\
=&\frac{b}{\Delta t}\left(q^{n},e_{x}^{n}-e_{x}^{n-1}\right)+\frac{b}{\Delta t}\left(e_{x}^{n},q^{n}-q^{n-1}\right)\\
&+b\left(q^{n}, \frac{(P_{h}^{0}u(t_{n}))_{x}-(P_{h}^{0}u(t_{n-1}))_{x}}{\Delta t}-(P_{h}^{0}u_{t}(t_{n}))_{x}\right)+b\left(e_{x}^{n}, \frac{P_{h}^{0}\phi(t_{n})-P_{h}^{0}\phi(t_{n-1})}{\Delta t}-P_{h}^{*}\phi_{t}(t_{n})\right)\\
=&\frac{b}{\Delta t}\left( (q^{n},e_{x}^{n})-(q^{n-1},e_{x}^{n-1})+(q^{n}-q^{n-1}, e_{x}^{n}-e_{x}^{n-1})\right)\\
&+b\left(q^{n}, \frac{(P_{h}^{0}u(t_{n}))_{x}-(P_{h}^{0}u(t_{n-1}))_{x}}{\Delta t}-(P_{h}^{0}u_{t}(t_{n}))_{x}\right)+b\left(e_{x}^{n}, \frac{P_{h}^{0}\phi(t_{n})-P_{h}^{0}\phi(t_{n-1})}{\Delta t}-P_{h}^{0}\phi_{t}(t_{n})\right),
\end{aligned}
\end{equation}
and 
\begin{equation}\label{3.7}
\begin{aligned}
&\mu\left(e_{x}^{n},\hat{e}_{x}^{n}\right)\\
=&\mu\left(e_{x}^{n},\frac{e_{x}^{n}-e_{x}^{n-1}}{\Delta t}+\frac{(P_{h}^{0}u(t_{n}))_{x}-(P_{h}^{0}u(t_{n-1}))_{x}}{\Delta t}-(P_{h}^{0}u_{t}(t_{n}))_{x}\right)\\
=&\frac{\mu}{2\Delta t}\left(||e_{x}^{n}-e_{x}^{n-1}||^{2}+||e_{x}^{n}||^{2}-||e_{x}^{n-1}||^{2}\right)+\mu\left(e_{x}^{n},\frac{(P_{h}^{0}u(t_{n}))_{x}-(P_{h}^{0}u(t_{n-1}))_{x}}{\Delta t}-(P_{h}^{0}u_{t}(t_{n}))_{x}\right),
\end{aligned}
\end{equation}
similarly 
\begin{equation}\label{3.8}
\begin{aligned}
&\delta\left(q_{x}^{n},\hat{q}_{x}^{n}\right)\\
=&\delta\left(q_{x}^{n},\frac{q_{x}^{n}-q_{x}^{n-1}}{\Delta t}+\frac{(P_{h}^{0}\phi(t_{n}))_{x}-(P_{h}^{0}\phi(t_{n-1}))_{x}}{\Delta t}-(P_{h}^{0}\phi_{t}(t_{n}))_{x}\right)\\
=&\frac{\delta}{2\Delta t}\left(||q_{x}^{n}-q_{x}^{n-1}||^{2}+||q_{x}^{n}||^{2}-||q_{x}^{n-1}||^{2}\right)+\delta\left(q_{x}^{n},\frac{(P_{h}^{0}\phi(t_{n}))_{x}-(P_{h}^{0}\phi(t_{n-1}))_{x}}{\Delta t}-(P_{h}^{0}\phi_{t}(t_{n}))_{x}\right),
\end{aligned}
\end{equation}
and 
\begin{equation}\label{3.9}
\begin{aligned}
&\xi\left(q^{n},\hat{q}^{n}\right)\\
=&\xi\left(q^{n},\frac{q^{n}-q^{n-1}}{\Delta t}+\frac{P_{h}^{0}\phi(t_{n})-P_{h}^{0}\phi(t_{n-1})}{\Delta t}-P_{h}^{0}\phi_{t}(t_{n})\right)\\
=&\frac{\xi}{2\Delta t}\left(||q^{n}-q^{n-1}||^{2}+||q^{n}||^{2}-||q^{n-1}||^{2}\right)+\xi\left(q^{n},\frac{P_{h}^{0}\phi(t_{n})-P_{h}^{0}\phi(t_{n-1})}{\Delta t}-P_{h}^{0}\phi_{t}(t_{n})\right).
\end{aligned}
\end{equation}
Now insert $\eqref{3.6}$, $\eqref{3.7}$, $\eqref{3.8}$ and $\eqref{3.9}$ in $\eqref{3.5}$ we get 
\begin{equation}\label{3.10}
\begin{aligned}
&\frac{\r}{2\Delta t}\left(||\hat{e}^{n}-\hat{e}^{n-1}||^{2}+||\hat{e}^{n}||^{2}-||\hat{e}^{n-1}||^{2}\right)+\frac{J\r}{2b\Delta t}\left(||\sigma\hat{e}^{n}||^{2}-||\sigma\hat{e}^{n-1}||^{2}+||\sigma(\hat{e}^{n}-\hat{e}^{n-1})||^{2}\right)\\
&+\frac{J\mu}{2b\Delta t}\left(||\hat{e}_{x}^{n}||^{2}-||\hat{e}_{x}^{n-1}||^{2}+||\hat{e}_{x}^{n}-\hat{e}_{x}^{n-1}||^{2}\right)+\frac{\mu}{2\Delta t}\left(||e_{x}^{n}-e_{x}^{n-1}||^{2}+||e_{x}^{n}||^{2}-||e_{x}^{n-1}||^{2}\right)\\
&\frac{\delta}{2\Delta t}\left(||q_{x}^{n}-q_{x}^{n-1}||^{2}+||q_{x}^{n}||^{2}-||q_{x}^{n-1}||^{2}\right)+\frac{\xi}{2\Delta t}\left(||q^{n}-q^{n-1}||^{2}+||q^{n}||^{2}-||q^{n-1}||^{2}\right)\\
&\frac{\a}{2\Delta t}\left(||R^{n}-R^{n-1}||^{2}+|| R^{n}||^{2}-|| R^{n-1}||^{2}\right)+\frac{b}{\Delta t}\left( (q^{n},e_{x}^{n})-(q^{n-1},e_{x}^{n-1})+(q^{n}-q^{n-1}, e_{x}^{n}-e_{x}^{n-1})\right)\\
&+\kappa||R_{x}^{n}||^{2}+k||R^{n}||^{2}\\
&=\r\left(u_{tt}(t_{n})-\frac{P_{h}^{0}u_{t}(t_{n})-P_{h}^{0}u_{t}(t_{n-1})}{\Delta t},\hat{e}^{n}\right)+\frac{J\r}{b}\left(u_{tt}(t_{n})-\frac{P_{h}^{0}u_{t}(t_{n})-P_{h}^{0}u_{t}(t_{n-1})}{\Delta t},\sigma^{2}\hat{e}^{n}\right)\\
&+\frac{J\mu}{b}\left(u_{ttx}(t_{n})-\frac{(P_{h}^{0}u_{t}(t_{n}))_{x}-(P_{h}^{0}u_{t}(t_{n-1}))_{x}}{\Delta t},\hat{e}_{x}^{n}\right)+\a\left(w_{t}(t_{n})-\frac{P_{h}^{0}w(t_{n})-P_{h}^{0}w(t_{n-1})}{\Delta t}, R^{n}\right)\\
&-\mu\left(e_{x}^{n},\frac{(P_{h}^{0}u(t_{n}))_{x}-(P_{h}^{0}u(t_{n-1}))_{x}}{\Delta t}-(P_{h}^{0}u_{t}(t_{n}))_{x}\right)-\delta\left(q_{x}^{n},\frac{(P_{h}^{0}\phi(t_{n}))_{x}-(P_{h}^{0}\phi(t_{n-1}))_{x}}{\Delta t}-(P_{h}^{0}\phi_{t}(t_{n}))_{x}\right)\\
&-\xi\left(q^{n},\frac{P_{h}^{0}\phi(t_{n})-P_{h}^{0}\phi(t_{n-1})}{\Delta t}-P_{h}^{0}\phi_{t}(t_{n})\right)-b\left(q^{n}, \frac{(P_{h}^{0}u(t_{n}))_{x}-(P_{h}^{0}u(t_{n-1}))_{x}}{\Delta t}-(P_{h}^{0}u_{t}(t_{n}))_{x}\right)\\
&-b\left(e_{x}^{n}, \frac{P_{h}^{0}\phi(t_{n})-P_{h}^{0}\phi(t_{n-1})}{\Delta t}-P_{h}^{0}\phi_{t}(t_{n})\right)\\
&+\mu\left(u_{x}(t_{n})-(P_{h}^{0}u(t_{n}))_{x},\hat{e}_{x}^{n}\right)+\delta\left(\phi_{x}(t_{n})-(P_{h}^{0}\phi(t_{n}))_{x},\hat{q}_{x}^{n}\right)+b\left(\phi(t_{n})-P_{h}^{0}\phi(t_{n}),\hat{e}_{x}^{n}\right)+b\left(u_{x}(t_{n})-(P_{h}^{0}u(t_{n}))_{x},\hat{q}^{n}\right)\\
&+\xi(\phi(t_{n})-P_{h}^{0}\phi(t_{n}),\hat{q}^{n})+d\left(w_{x}(t_{n})-(P_{h}^{0}w(t_{n}))_{x},\hat{q}^{n}\right)+\kappa\left(w_{x}(t_{n})-(P_{h}^{0}w(t_{n}))_{x}, R_{x}^{n}\right)\\
&+d\left(\phi_{tx}(t_{n})-(P_{h}^{0}\phi_{t}(t_{n}))_{x}, R^{n}\right)+k(w-P_{h}^{0}w(t_{n}),R^{n}).
\end{aligned}
\end{equation}

\textbf{Step 5.} Since $\mu\xi-b^{2}>0$ we get 
$$\frac{\mu}{2\Delta t}||e_{x}^{n}-e_{x}^{n-1}||^{2}+\frac{\xi}{2\Delta t}||q^{n}-q^{n-1}||^{2}+\frac{b}{\Delta t}(q^{n}-q^{n-1}, e_{x}^{n}-e_{x}^{n-1})>0,$$
and by taking into consideration the positivity of the terms $||\hat{e}^{n}-\hat{e}^{n-1}||^{2}$, $||\hat{e}_{x}^{n}-\hat{e}_{x}^{n-1}||^{2}$, $||(\hat{e}^{n}-\hat{e}^{n-1})_{t}||^{2}$, $||q_{x}^{n}-q_{x}^{n-1}||^{2}$ and $||R^{n}-R^{n-1}||^{2}$, we arrive at the following inequality
 \begin{equation}\label{3.11}
\begin{aligned}
&\frac{\r}{2\Delta t}\left(||\hat{e}^{n}||^{2}-||\hat{e}^{n-1}||^{2}\right)+\frac{J\r}{2b\Delta t}\left(||\sigma\hat{e}^{n}||^{2}-||\sigma\hat{e}^{n-1}||^{2}\right)\\
&+\frac{J\mu}{2b\Delta t}\left(||\hat{e}_{x}^{n}||^{2}-||\hat{e}_{x}^{n-1}||^{2}\right)+\frac{\mu}{2\Delta t}\left(||e_{x}^{n}||^{2}-||e_{x}^{n-1}||^{2}\right)\\
&\frac{\delta}{2\Delta t}\left(||q_{x}^{n}||^{2}-||q_{x}^{n-1}||^{2}\right)+\frac{\xi}{2\Delta t}\left(||q^{n}||^{2}-||q^{n-1}||^{2}\right)\\
&\frac{\a}{2\Delta t}\left(|| R^{n}||^{2}-|| R^{n-1}||^{2}\right)+\frac{b}{\Delta t}\left( (q^{n},e_{x}^{n})-(q^{n-1},e_{x}^{n-1})\right)\\
&+\kappa||R_{x}^{n}||^{2}+k||R^{n}||^{2}\\
&<\r\left(u_{tt}(t_{n})-\frac{P_{h}^{0}u_{t}(t_{n})-P_{h}^{0}u_{t}(t_{n-1})}{\Delta t},\hat{e}^{n}\right)+\frac{J\r}{b}\left(u_{tt}(t_{n})-\frac{P_{h}^{0}u_{t}(t_{n})-P_{h}^{0}u_{t}(t_{n-1})}{\Delta t},\sigma^{2}\hat{e}^{n}\right)\\
&+\frac{J\mu}{b}\left(u_{ttx}(t_{n})-\frac{(P_{h}^{0}u_{t}(t_{n}))_{x}-(P_{h}^{0}u_{t}(t_{n-1}))_{x}}{\Delta t},\hat{e}_{x}^{n}\right)+\a\left(w_{t}(t_{n})-\frac{P_{h}^{0}w(t_{n})-P_{h}^{0}w(t_{n-1})}{\Delta t}, R^{n}\right)\\
&-\mu\left(e_{x}^{n},\frac{(P_{h}^{0}u(t_{n}))_{x}-(P_{h}^{0}u(t_{n-1}))_{x}}{\Delta t}-(P_{h}^{0}u_{t}(t_{n}))_{x}\right)-\delta\left(q_{x}^{n},\frac{(P_{h}^{0}\phi(t_{n}))_{x}-(P_{h}^{0}\phi(t_{n-1}))_{x}}{\Delta t}-(P_{h}^{0}\phi_{t}(t_{n}))_{x}\right)\\
&-\xi\left(q^{n},\frac{P_{h}^{0}\phi(t_{n})-P_{h}^{0}\phi(t_{n-1})}{\Delta t}-P_{h}^{0}\phi_{t}(t_{n})\right)-b\left(q^{n}, \frac{(P_{h}^{0}u(t_{n}))_{x}-(P_{h}^{0}u(t_{n-1}))_{x}}{\Delta t}-(P_{h}^{0}u_{t}(t_{n}))_{x}\right)\\
&-b\left(e_{x}^{n}, \frac{P_{h}^{0}\phi(t_{n})-P_{h}^{0}\phi(t_{n-1})}{\Delta t}-P_{h}^{0}\phi_{t}(t_{n})\right)\\
&+\mu\left(u_{x}(t_{n})-(P_{h}^{0}u(t_{n}))_{x},\hat{e}_{x}^{n}\right)+\delta\left(\phi_{x}(t_{n})-(P_{h}^{0}\phi(t_{n}))_{x},\hat{q}_{x}^{n}\right)+b\left(\phi(t_{n})-P_{h}^{0}\phi(t_{n}),\hat{e}_{x}^{n}\right)+b\left(u_{x}(t_{n})-(P_{h}^{0}u(t_{n}))_{x},\hat{q}^{n}\right)\\
&+\xi(\phi(t_{n})-P_{h}^{0}\phi(t_{n}),\hat{q}^{n})+d\left(w_{x}(t_{n})-(P_{h}^{0}w(t_{n}))_{x},\hat{q}^{n}\right)+\kappa\left(w_{x}(t_{n})-(P_{h}^{0}w(t_{n}))_{x}, R_{x}^{n}\right)\\
&+d\left(\phi_{tx}(t_{n})-(P_{h}^{0}\phi_{t}(t_{n}))_{x}, R^{n}\right)+k(w-P_{h}^{0}w(t_{n}),R^{n}).
\end{aligned}
\end{equation}
Now we want to eliminate the damping parameters from the left side of inequality $\eqref{3.11}$,  denote by 
$$\Theta_{1}=k(w-P_{h}^{0}w(t_{n}),R^{n})$$
By using Young's inequality for all $\eps_{5}>0$  
\begin{equation*}
\Theta_{1}\leq \frac{k}{2\eps_{5}}||w-P_{h}^{0}w(t_{n})||^{2}+\frac{k\eps_{5}}{2}||R^{n}||^{2}
\end{equation*}
Choose $\eps_{5}=2$ then 
\begin{equation}\label{3.12}
\Theta_{1}\leq \frac{k}{4}||w-P_{h}^{0}w(t_{n})||^{2}+k||R^{n}||^{2}.
\end{equation}
 Denote by 
$$\Theta_{2}=\kappa(w_{x}-(P_{h}^{0}w(t_{n}))_{x},R_{x}^{n})$$
By using Young's inequality for all $\eps_{6}>0$  
\begin{equation*}
\Theta_{2}\leq \frac{\kappa}{2\eps_{6}}||w_{x}-(P_{h}^{0}w(t_{n}))_{x}||^{2}+\frac{\kappa\eps_{6}}{2}||R_{x}^{n}||^{2}
\end{equation*}
Choose $\eps_{6}=2$ then 
\begin{equation}\label{3.13}
\Theta_{2}\leq \frac{\kappa}{4}||w_{x}-(P_{h}^{0}w(t_{n}))_{x}||^{2}+\kappa||R_{x}^{n}||^{2}.
\end{equation}
Insert inequalities $\eqref{3.12}$ and $\eqref{3.13}$ in $\eqref{3.11}$ we obtain 
\begin{equation}\label{3.11*}
\begin{aligned}
&\frac{\r}{2\Delta t}\left(||\hat{e}^{n}||^{2}-||\hat{e}^{n-1}||^{2}\right)+\frac{J\r}{2b\Delta t}\left(||\sigma\hat{e}^{n}||^{2}-||\sigma\hat{e}^{n-1}||^{2}\right)\\
&+\frac{J\mu}{2b\Delta t}\left(||\hat{e}_{x}^{n}||^{2}-||\hat{e}_{x}^{n-1}||^{2}\right)+\frac{\mu}{2\Delta t}\left(||e_{x}^{n}||^{2}-||e_{x}^{n-1}||^{2}\right)\\
&\frac{\delta}{2\Delta t}\left(||q_{x}^{n}||^{2}-||q_{x}^{n-1}||^{2}\right)+\frac{\xi}{2\Delta t}\left(||q^{n}||^{2}-||q^{n-1}||^{2}\right)\\
&\frac{\a}{2\Delta t}\left(|| R^{n}||^{2}-|| R^{n-1}||^{2}\right)+\frac{b}{\Delta t}\left( (q^{n},e_{x}^{n})-(q^{n-1},e_{x}^{n-1})\right)\\
&<\r\left(u_{tt}(t_{n})-\frac{P_{h}^{0}u_{t}(t_{n})-P_{h}^{0}u_{t}(t_{n-1})}{\Delta t},\hat{e}^{n}\right)+\frac{J\r}{b}\left(u_{tt}(t_{n})-\frac{P_{h}^{0}u_{t}(t_{n})-P_{h}^{0}u_{t}(t_{n-1})}{\Delta t},\sigma^{2}\hat{e}^{n}\right)\\
&+\frac{J\mu}{b}\left(u_{ttx}(t_{n})-\frac{(P_{h}^{0}u_{t}(t_{n}))_{x}-(P_{h}^{0}u_{t}(t_{n-1}))_{x}}{\Delta t},\hat{e}_{x}^{n}\right)+\a\left(w_{t}(t_{n})-\frac{P_{h}^{0}w(t_{n})-P_{h}^{0}w(t_{n-1})}{\Delta t}, R^{n}\right)\\
&+\mu\left(e_{x}^{n},(P_{h}^{0}u_{t}(t_{n}))_{x}-\frac{(P_{h}^{0}u(t_{n}))_{x}-(P_{h}^{0}u(t_{n-1}))_{x}}{\Delta t}\right)+\delta\left(q_{x}^{n},(P_{h}^{0}\phi_{t}(t_{n}))_{x}-\frac{(P_{h}^{0}\phi(t_{n}))_{x}-(P_{h}^{0}\phi(t_{n-1}))_{x}}{\Delta t}\right)\\
&+\xi\left(q^{n},P_{h}^{0}\phi_{t}(t_{n})-\frac{P_{h}^{0}\phi(t_{n})-P_{h}^{0}\phi(t_{n-1})}{\Delta t}\right)+b\left(q^{n}, (P_{h}^{0}u_{t}(t_{n}))_{x}-\frac{(P_{h}^{0}u(t_{n}))_{x}-(P_{h}^{0}u(t_{n-1}))_{x}}{\Delta t}\right)\\
&+b\left(e_{x}^{n},P_{h}^{0}\phi_{t}(t_{n})- \frac{P_{h}^{0}\phi(t_{n})-P_{h}^{0}\phi(t_{n-1})}{\Delta t}\right)\\
&+\mu\left(u_{x}(t_{n})-(P_{h}^{0}u(t_{n}))_{x},\hat{e}_{x}^{n}\right)+\delta\left(\phi_{x}(t_{n})-(P_{h}^{0}\phi(t_{n}))_{x},\hat{q}_{x}^{n}\right)+b\left(\phi(t_{n})-P_{h}^{0}\phi(t_{n}),\hat{e}_{x}^{n}\right)+b\left(u_{x}(t_{n})-(P_{h}^{0}u(t_{n}))_{x},\hat{q}^{n}\right)\\
&+\xi(\phi(t_{n})-P_{h}^{0}\phi(t_{n}),\hat{q}^{n})+d\left(w_{x}(t_{n})-(P_{h}^{0}w(t_{n}))_{x},\hat{q}^{n}\right)+\frac{\kappa}{4}||w_{x}-(P_{h}^{0}w(t_{n}))_{x}||^{2}\\
&+d\left(\phi_{tx}(t_{n})-(P_{h}^{0}\phi_{t}(t_{n}))_{x}, R^{n}\right)+\frac{k}{4}||w-P_{h}^{0}w(t_{n})||^{2}.
\end{aligned}
\end{equation}
\textbf{Step 7.} Now let 
$$Z_{n}=\r||\hat{e}^{n}||^{2}+\frac{J\r}{b}||\sigma\hat{e}^{n}||^{2}+\frac{J\mu}{b}||\hat{e}_{x}^{n}||^{2}+\mu||e_{x}^{n}||^{2}+\delta||q_{x}^{n}||^{2}+\xi||q^{n}||^{2}+\a||R^{n}||^{2}.$$
By Young's inequality there exists a constant $c$ such that 
\begin{equation}\label{3.15}
Z_{n}-Z_{n-1}+b (q^{n},e_{x}^{n})-b(q^{n-1},e_{x}^{n-1})\leq 2c\Delta t(Z_{n}+K_{n}),
\end{equation}
where 
\begin{equation*}
\begin{aligned}
K_{n}&=\underbrace{\left\|u_{tt}(t_{n})-\frac{P_{h}^{0}u_{t}(t_{n})-P_{h}^{0}u_{t}(t_{n-1})}{\Delta t}\right\|^{2}}_{ I_{1}}+\underbrace{\left\|u_{tt}(t_{n})-\frac{P_{h}^{0}u_{t}(t_{n})-P_{h}^{0}u_{t}(t_{n-1})}{\Delta t}\right\|^{2}}_{I_{2}}\\
&+\underbrace{\left\|u_{ttx}(t_{n})-\frac{(P_{h}^{0}u_{t}(t_{n}))_{x}-(P_{h}^{0}u_{t}(t_{n-1}))_{x}}{\Delta t}\right\|^{2}}_{I_{3}}+\underbrace{\left\|(P_{h}^{0}u_{t}(t_{n}))_{x}-\frac{(P_{h}^{0}u(t_{n}))_{x}-(P_{h}^{0}u(t_{n-1}))_{x}}{\Delta t}\right\|^{2}}_{ I_{4}}\\
&+\underbrace{\left\|(P_{h}^{0}\phi_{t}(t_{n}))_{x}-\frac{(P_{h}^{0}\phi(t_{n}))_{x}-(P_{h}^{0}\phi(t_{n-1}))_{x}}{\Delta t}\right\|^{2}}_{I_{5}}+\underbrace{\left\|P_{h}^{0}\phi_{t}(t_{n})-\frac{P_{h}^{0}\phi(t_{n})-P_{h}^{0}\phi(t_{n-1})}{\Delta t}\right\|^{2}}_{ I_{6}}\\
&+\underbrace{\left\|(P_{h}^{0}u_{t}(t_{n}))_{x}-\frac{(P_{h}^{0}u(t_{n}))_{x}-(P_{h}^{0}u(t_{n-1}))_{x}}{\Delta t}\right\|^{2}}_{ I_{7}}+\underbrace{\left\|P_{h}^{0}\phi_{t}(t_{n})-\frac{P_{h}^{0}\phi(t_{n})-P_{h}^{0}\phi(t_{n-1})}{\Delta t}\right\|^{2}}_{ I_{8}}\\
&\underbrace{\left\|w_{t}(t_{n})-\frac{P_{h}^{0}w(t_{n})-P_{h}^{0}w(t_{n-1})}{\Delta t}\right\|^{2}}_{I_{9}}+\underbrace{\left\|u_{x}(t_{n})-(P_{h}^{0}u(t_{n}))_{x}\right\|^{2}}_{ I_{10}}+\underbrace{\left\|\phi_{x}(t_{n})-(P_{h}^{0}\phi(t_{n}))_{x}\right\|^{2}}_ {I_{11}}\\
&+\underbrace{\left\|\phi(t_{n})-P_{h}^{0}\phi(t_{n})\right\|^{2}}_ {I_{12}}+\underbrace{\left\|u_{x}(t_{n})-(P_{h}^{0}u(t_{n}))_{x}\right\|^{2}}_{I_{13}}+\underbrace{\left\|\phi(t_{n})-P_{h}^{0}\phi(t_{n})\right\|^{2}}_{ I_{14}}+\underbrace{\left\|\phi_{tx}(t_{n})-(P_{h}^{0}\phi_{t}(t_{n}))_{x}\right\|^{2}}_ {I_{15}}\\
&+\underbrace{\|w_{x}(t_{n})-(P_{h}^{0}w(t_{n}))_{x}\|^{2}}_{I_{16}}+\underbrace{\|w(t_{n})-P_{h}^{0}w(t_{n})\|^{2}}_{I_{17}}
\end{aligned}
\end{equation*}
Now by extending $\phi(t_{n-1})$ using Taylor formula near $t_{n}$ and using the linearity and properties of $P_{h}^{0}$ we get :
\begin{equation*}
\begin{aligned}
I_{6}&=\left\|P_{h}^{0}\phi_{t}(t_{n})-\frac{P_{h}^{0}\phi(t_{n})-P_{h}^{0}\phi(t_{n-1})}{\Delta t}\right\|^{2}\\
&\leq 2\left\|P_{h}^{0}\phi_{t}(t_{n})-\phi_{t}(t_{n})\right\|^{2}+2\left\|\phi_{t}(t_{n})-\frac{P_{h}^{0}\phi(t_{n})-P_{h}^{0}\phi(t_{n-1})}{\Delta t}\right\|^{2}\\
&\leq Ch^{2}\|\phi_{tx}(t_{n})\|^{2}+2\left\|\phi_{t}(t_{n})-\frac{P_{h}^{0}\phi(t_{n})-P_{h}^{0}\phi(t_{n-1})}{\Delta t}\right\|^{2}\\
&\leq Ch^{2}\|\phi_{tx}(t_{n})\|^{2}+2\left\|\phi_{t}(t_{n})-\frac{P_{h}^{0}\phi(t_{n})-P_{h}^{0}\phi(t_{n}-\Delta t)}{\Delta t}\right\|^{2}\\
&\leq Ch^{2}\|\phi_{tx}(t_{n})\|^{2}+2\left\|\phi_{t}(t_{n})-\frac{P_{h}^{0}\phi(t_{n})-P_{h}^{0}[\phi(t_{n})-\Delta t\phi_{t}(t_{n})+\frac{\Delta t^{2}}{2}\phi_{tt}(\xi_{n})]}{\Delta t}\right\|^{2}\\
&\leq Ch^{2}\|\phi_{tx}(t_{n})\|^{2}+2\left\|\phi_{t}(t_{n})-P_{h}^{0}\phi_{t}(t_{n})-\frac{\Delta t^{2}}{2}\phi_{tt}(\xi_{n})\right\|^{2}\\
&\leq  Ch^{2}\|\phi_{tx}(t_{n})\|^{2}+ c\Delta t^{2}\left\|\phi_{tt}(\xi_{n})\right\|.
\end{aligned}
\end{equation*}
We introduce 
$$f(\xi_{n})=\left\|\phi_{tt}(\xi_{n})\right\|. $$
Then by Galiardo-Neirenberg inequality, 
\begin{equation*}
\begin{aligned}
f(\xi_{n})&\leq C \|f\|_{L^{2}(0,T)}\|f_{t}\|_{L^{2}(0,T)}\\
&\leq C\left( \int_{0}^{T}\left\|\phi_{tt}\right\|^{2}dt\right)^{\frac{1}{2}}\left( \int_{0}^{T}\left\|\phi_{ttt}\right\|^{2}dt\right)^{\frac{1}{2}}\\
&\leq C\left\|\phi\right\|_{H^{2}(0,T;L^{2})}\left\|\phi\right\|_{H^{3}(0,T;L^{2})}.
\end{aligned}
\end{equation*}
Similarly we show that 
$$\|\phi_{tx}(t_{n})\|\leq C\left\|\phi\right\|_{H^{1}(0,T;H^{1})}\left\|\phi\right\|_{H^{2}(0,T;H^{1})}.$$
Therefore 
$$I_{6}\leq C(h^{2}+\Delta t^{2})\left\|\phi\right\|_{H^{3}(0,T;H^{1})}^{2}$$
With the same manner 
\begin{equation*}
\begin{aligned}
&I_{1}\leq C(h^{2}+\Delta t^{2})\left\|u\right\|_{H^{4}(0,T;H^{1})}^{2}\\
&I_{2}\leq C(h^{2}+\Delta t^{2})\left\|u\right\|_{H^{4}(0,T;H^{1})}^{2}\\
&I_{3}\leq C(h^{2}+\Delta t^{2})\left\|u\right\|_{H^{4}(0,T;H^{2})}^{2}\\
&I_{4}\leq C(h^{2}+\Delta t^{2})\left\|u\right\|_{H^{3}(0,T;H^{2})}^{2}\\
&I_{5}\leq C(h^{2}+\Delta t^{2})\left\|\phi\right\|_{H^{3}(0,T;H^{2})}^{2}\\
&I_{7}\leq C(h^{2}+\Delta t^{2})\left\|u\right\|_{H^{3}(0,T;H^{2})}^{2}\\
&I_{8}\leq C(h^{2}+\Delta t^{2})\left\|\phi\right\|_{H^{3}(0,T;H^{1})}^{2}\\
&I_{9}\leq C(h^{2}+\Delta t^{2})\left\|w\right\|_{H^{3}(0,T;H^{1})}^{2}\\
&I_{10}\leq Ch^{2}\left\|u\right\|_{L^{2}(0,T;H^{2})}^{2}\\
&I_{11}\leq Ch^{2}\left\|\phi\right\|_{L^{2}(0,T;H^{2})}^{2}\\
&I_{12}\leq Ch^{2}\left\|\phi\right\|_{L^{2}(0,T;H^{1})}^{2}\\
&I_{13}\leq Ch^{2}\left\|u\right\|_{L^{2}(0,T;H^{2})}^{2}\\
&I_{14}\leq Ch^{2}\left\|\phi\right\|_{L^{2}(0,T;H^{1})}^{2}\\
&I_{15}\leq Ch^{2}\left\|\phi\right\|_{H^{1}(0,T;H^{2})}^{2}\\
&I_{16}\leq Ch^{2}\left\|w\right\|_{L^{2}(0,T;H^{2})}^{2}\\
&I_{17}\leq Ch^{2}\left\|w\right\|_{L^{2}(0,T;H^{1})}^{2}\\
\end{aligned}
\end{equation*}
Now sum equation $\eqref{3.15}$ over $n$ we get 
$$Z_{j}-Z_{0}+b (q^{j},e_{x}^{j})-b(q^{0},e_{x}^{0})\leq 2c\Delta t\sum_{n=0}^{j}(Z_{n}+K_{n}).$$
Since $Z_{0}=0$, we have 
$$Z_{j}< 2c\Delta t\sum_{n=0}^{j}(Z_{n}+K_{n})-b (q^{j},e_{x}^{j}).$$
By applying Young's inequality there is $\nu$ such that 
$$Z_{j}< 2c\Delta t\sum_{n=0}^{j}(Z_{n}+K_{n})+\frac{b\nu}{2}\|q^{j}\|^{2}+\frac{b}{2\nu}\|e_{x}^{j}\|^{2},$$
where $\nu$ is chosen in the following way 
$$\frac{b}{2\mu}<\nu<\frac{2\xi}{b},$$
hence 
$$CZ_{j}<2c\Delta t\sum_{n=0}^{j}Z_{n}+c\Delta t(\Delta t^{2}+h^{2}).$$
Finally apply Grownwall's inequality the proof is completed $\hfill{\blacksquare}$
\section{Numerical Simulations}

Knowing that $\hat{u}_{h}^{n}=\frac{u_{h}^{n}-u_{h}^{n-1}}{\Delta t}$, $\hat{%
\phi}_{h}^{n}=\frac{\phi _{h}^{n}-\phi _{h}^{n-1}}{\Delta t}$ and $\hat{w}%
_{h}^{n}=\frac{w_{h}^{n}-w_{h}^{n-1}}{\Delta t}$ and plugging them in system \ref{NP}
we get,

\[
\left\{ 
\begin{array}{l}
\frac{\rho }{\left( \Delta t\right) ^{2}}\left(
u_{h}^{n}-2u_{h}^{n-1}+u_{h}^{n-2},\bar{u}_{h}\right) +\mu \left( u_{hx}^{n},%
\bar{u}_{hx}\right) +b\left( \phi _{h}^{n},\bar{u}_{hx}\right) =0, \\ 
\frac{J}{\left( \Delta t\right) ^{2}}\left( \phi _{h}^{n}-2\phi
_{h}^{n-1}+\phi _{h}^{n-2},\bar{\phi}_{hx}\right) +\delta \left( \phi
_{hx}^{n},\bar{\phi}_{hx}\right) +b\left( u_{hx}^{n},\bar{\phi}_{h}\right)
+\xi \left( \phi _{h}^{n},\bar{\phi}_{h}\right) +d\left( w_{hx}^{n},\bar{\phi%
}_{h}\right) =0 \\ 
\frac{\alpha }{\Delta t}\left( w_{h}^{n}-w_{h}^{n-1},\bar{w}_{h}\right)
+\kappa \left( w_{hx}^{n},\bar{w}_{hx}\right) +\frac{d}{\Delta t}\left( \phi
_{hx}^{n}-\phi _{hx}^{n-1},\bar{w}_{h}\right) +k\left( w_{h}^{n},\bar{w}%
_{h}\right) =0%
\end{array}%
\right. 
\]

\noindent note that by the finite element theory, $u_{h}^{n}=\sum_{i=1}^{s}a_{i}^{n}%
\psi _{i}$, $\phi _{h}^{n}=\sum_{i=1}^{s}b_{i}^{n}\psi _{i}$ and $%
w_{h}^{n}=\sum_{i=1}^{s}c_{i}^{n}\psi _{i}$ where $\psi _{i}$ are bases of
the finite space $S_{h}^{0}$. Taking $\bar{u}_{h}=\psi _{j}$ we get

\[
\left\{ 
\begin{array}{l}\label{matrices}
\frac{\rho }{\left( \Delta t\right) ^{2}}A^{n}Z+\mu A^{n}T+bB^{n}Y=\frac{%
2\rho }{\left( \Delta t\right) ^{2}}A^{n-1}Z-\frac{\rho }{\left( \Delta
t\right) ^{2}}A^{n-2}Z, \\ 
\frac{J}{\left( \Delta t\right) ^{2}}A^{n}Y+\delta B^{n}T+bA^{n}X+\xi
B^{n}Z+dC^{n}X=\frac{2J}{\left( \Delta t\right) ^{2}}A^{n-1}Y-\frac{J}{%
\left( \Delta t\right) ^{2}}A^{n-2}Y \\ 
\frac{\alpha }{\Delta t}C^{n}Z+\kappa C^{n}T+\frac{d}{\Delta t}%
B^{n}X+kC^{n}Z=\frac{\alpha }{\Delta t}C^{n-1}Z+\frac{d}{\Delta t}B^{n-1}X%
\end{array}%
\right. 
\]

\noindent Where the vectors $A^{n}$, $B^{n}$ and $C^{n}$ are given by 
\begin{eqnarray*}
A^{n} &=&\left\{ a_{i}^{n}\right\} \\
B^{n} &=&\left\{ b_{i}^{n}\right\} \\
C^{n} &=&\left\{ c_{i}^{n}\right\}
\end{eqnarray*}

\noindent and the matrices $Z$, $X$, $Y$, and $T$ are given by

\begin{eqnarray*}
X &=&\left( \psi _{ix},\psi _{j}\right) \\
Y &=&\left( \psi _{i},\psi _{jx}\right) \\
Z &=&\left( \psi _{i},\psi _{j}\right) \\
T &=&\left( \psi _{ix},\psi _{jx}\right)
\end{eqnarray*}
We solve \ref{matrices} using the following initial and
physical data, $\rho=d=\alpha=b=\xi=j=0.001$, $k=1$ and $\mu=0.01$. The space discritization $\Delta x=\frac{1}{11}$ and the time dicritization $\Delta t=\frac{1}{22}$ with total time $T=25$. The initial data $u_{h}^{0}=\phi_{h}^{0}=w_{h}^{0}=u_{h}^{1}=\phi_{h}^{1}=(1-x)x$
\begin{figure}
	\centering
	\includegraphics[width=0.50\textwidth]{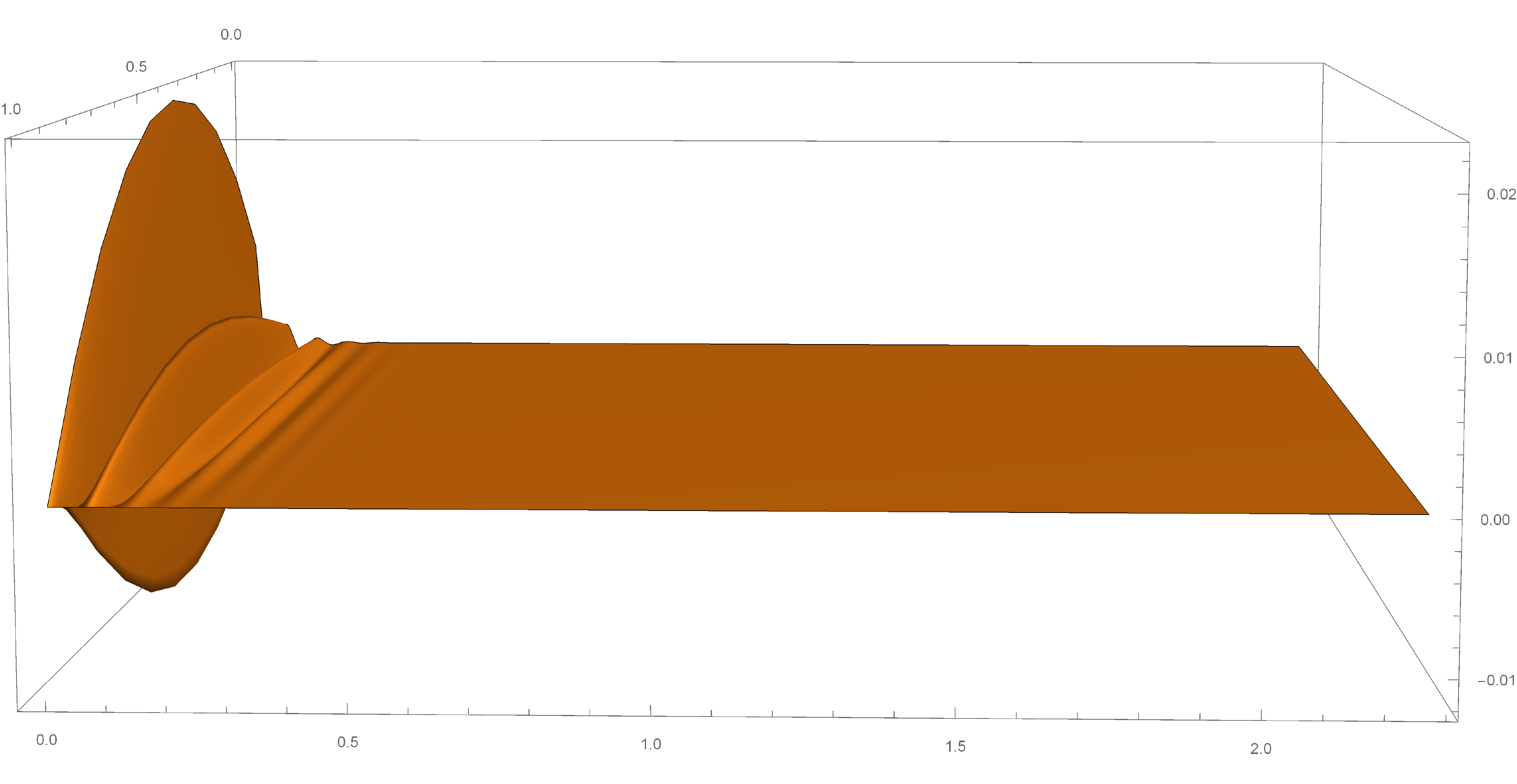}
	\caption{$u $ as function of $x$ and $t$.}
	\label{xip37}
\end{figure}
\begin{figure}
	\centering
	\includegraphics[width=0.50\textwidth]{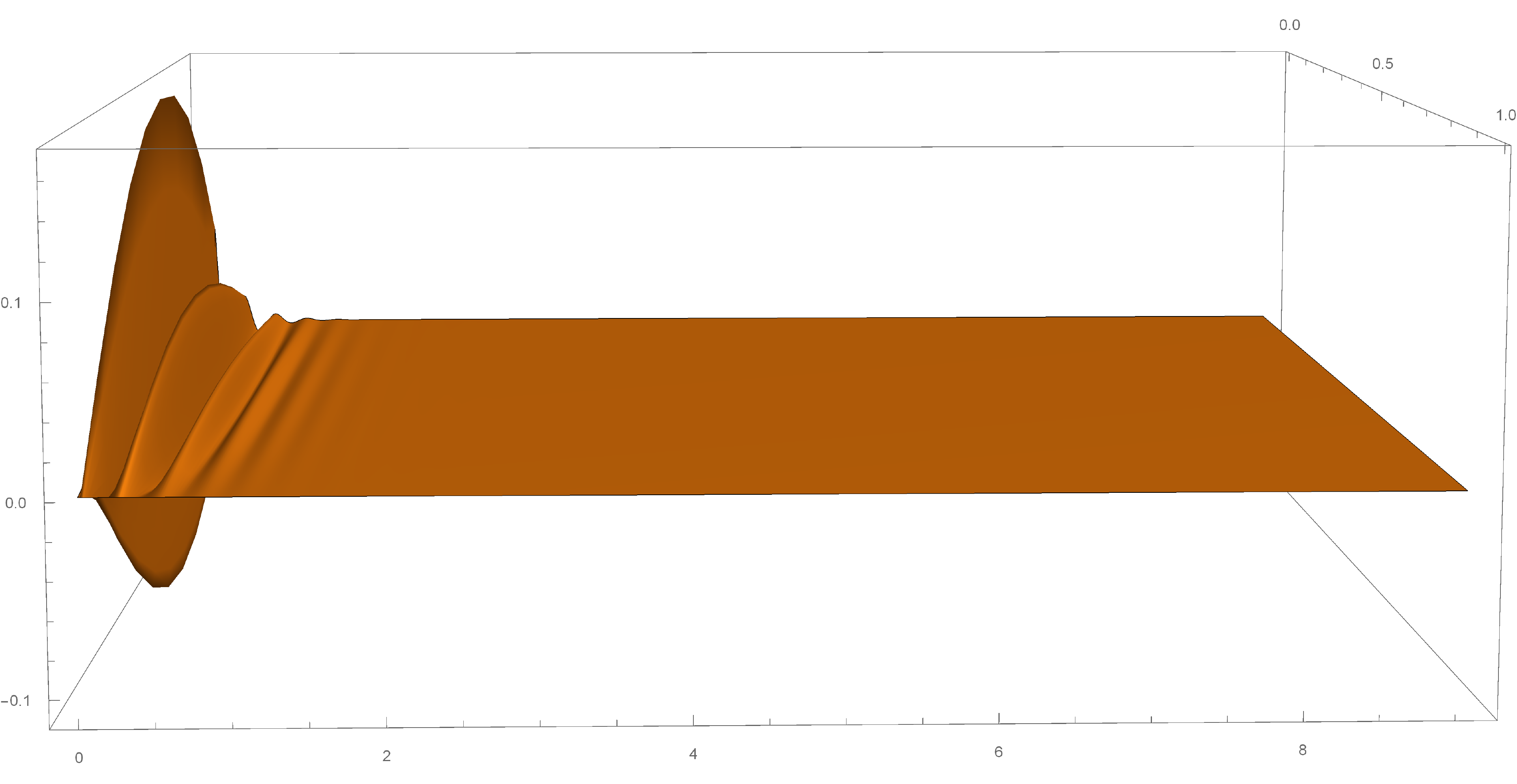}
	\caption{$\phi $ as function of $x$ and $t$.}
	\label{xip38}
\end{figure}	
\begin{figure}
	\centering
	\includegraphics[width=0.50\textwidth]{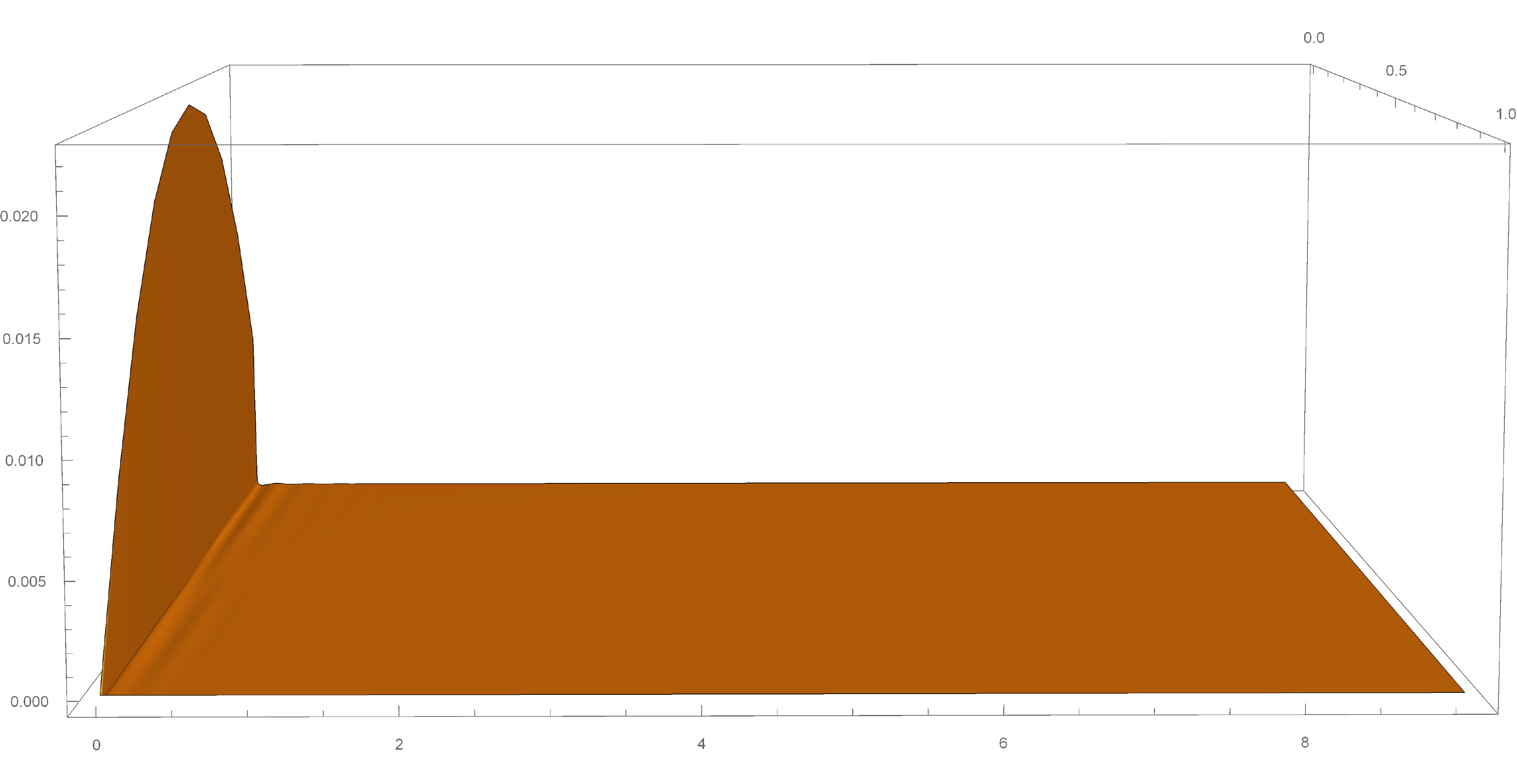}
	\caption{$w $ as function of $x$ and $t$.}
	\label{xip39}
\end{figure}
\begin{figure}
	\centering
	\includegraphics[width=0.50\textwidth]{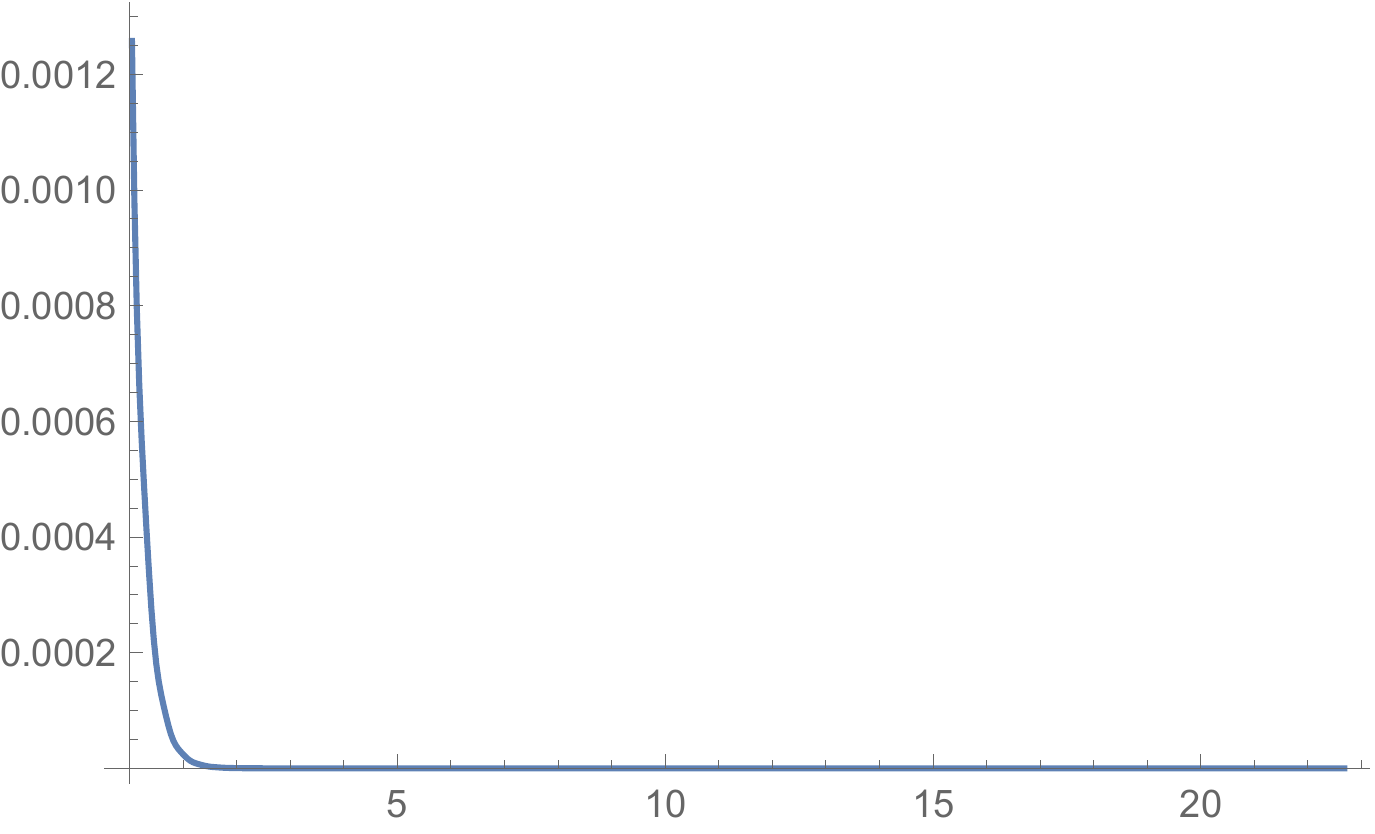}
	\caption{Energy as function of time.}
	\label{xip311}
\end{figure}
\begin{figure}
	\centering
	\includegraphics[width=0.50\textwidth]{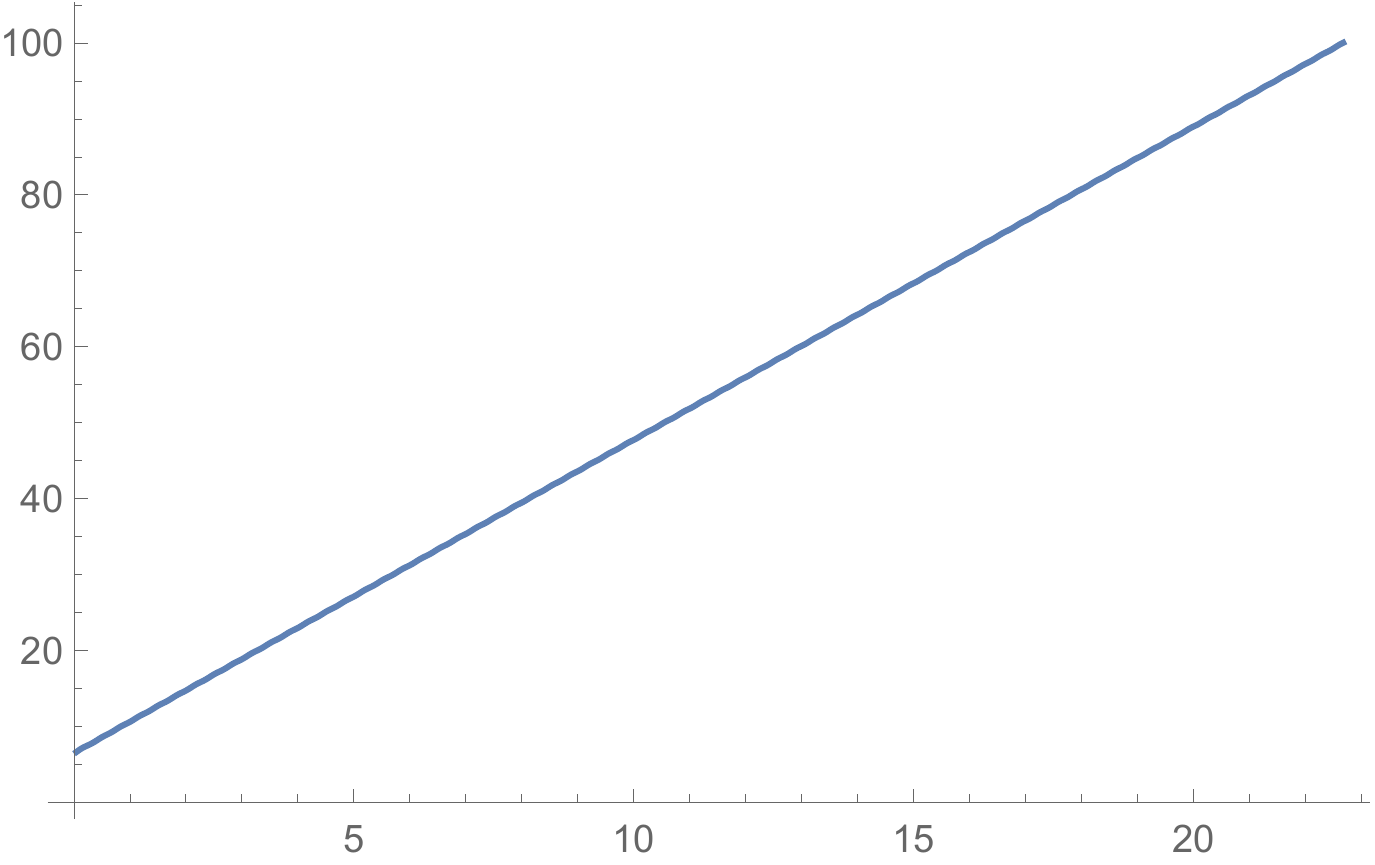}
	\caption{ $-log[energy] $ as function of time.}
	\label{xip312}
\end{figure}

\end{document}